\documentstyle[fleqn]{article}
\newif\ifams

\amstrue		%amssym.def and amssym.tex are available
\ifams
   \input{amssym.def}
   \input{amssym.tex}
\fi

\newcommand{\myfootnote}[1]{\footnote{#1}{\rm)}}

\newcommand{\proved}{\hspace*{\fill}$\Box$}
\newcommand{\PP}{I\!\!P}
\newcommand{\PPV}{\check{I\!\!P}}

\ifams
   \newcommand{\complex}{\Bbb C}
\else
   \def\complex{{\mathchoice {\setbox0=\hbox{$\displaystyle\rm C$}\hbox{\hbox
      to0pt{\kern0.4\wd0\vrule height0.9\ht0\hss}\box0}}
      {\setbox0=\hbox{$\textstyle\rm C$}\hbox{\hbox
      to0pt{\kern0.4\wd0\vrule height0.9\ht0\hss}\box0}}
      {\setbox0=\hbox{$\scriptstyle\rm C$}\hbox{\hbox
      to0pt{\kern0.4\wd0\vrule height0.9\ht0\hss}\box0}}
      {\setbox0=\hbox{$\scriptscriptstyle\rm C$}\hbox{\hbox
      to0pt{\kern0.4\wd0\vrule height0.9\ht0\hss}\box0}}}}
\fi

\newcommand{\lhookrightarrow}{\lhook\joinrel\longrightarrow}
\newcommand{\ratarrow}{\relbar\relbar\rightarrow}

\ifams
\else
   \newcommand{\ltimes}{\:{\vrule height 1ex width 0.25pt depth -0.1ex%
      \mkern -2.5mu\times}}
\fi

\newcommand{\spec}{\mbox{Spec}}

\newcommand{\ti}[1]{\widetilde{#1}}
\renewcommand{\bar}[1]{\overline{#1}}
\newcommand{\Grass}{\mbox{Grass}}
\newcommand{\Fano}{\mbox{Fano}}
\newcommand{\Bisec}{\mbox{Bisec}}
\newcommand{\Pic}{\mbox{Pic}}
\renewcommand{\emptyset}{\big/ \!\!\!\!\!\!\:\bigcirc}
\newcommand{\symm}{\stackrel{\mbox{\tiny sym}}{\times}}
\newcommand{\hateq}{\,\widehat{=}\,}
\newcommand{\POne}{\mbox{$P_{\!1}$}}

\newcommand{\OA}{\mbox{$\cal A$}}
\newcommand{\OB}{\mbox{$\cal B$}}
\newcommand{\OC}{\mbox{$\cal C$}}

\newcommand{\Esechs}{\mbox{${\bf E}_6$}}
\newcommand{\OO}{{\cal O}}
\newcommand{\R}{\mbox{$I\!\!R$}}
\newcommand{\YO}{\mbox{$Y_0$}}
\newcommand{\YF}{\mbox{$Y_F$}}

\newcommand{\YY}{\mbox{$Y_F^2$}}
\newcommand{\X}{\mbox{$X$}}

\ifams
   \newcommand{\Z}{\mbox{$\Bbb Z$}}
\else
   \newcommand{\Z}{\mbox{$Z\!\!\!Z$}}
\fi

\newenvironment{morph}{\[\renewcommand{\arraystretch}{1.5}\begin{array}{rccc}}%
{\end{array}\]}
\newenvironment{proof}{{\bf Proof:}}{\proved}

\begin{document}

\newtheorem{satz}{Proposition}[section]
\newtheorem{lemma}[satz]{Lemma}
\newtheorem{folg}[satz]{Corollary}
\newtheorem{theorem}[satz]{Theorem}

\title{Conics Touching a Quartic Surface with 13 Nodes}
\author{Ingo Hadan}
\date{February 21, 1996}
\maketitle

\begin{abstract}
\noindent
For a given real quartic surface in complex $\PP ^3$ that has exactly 13
ordinary nodes as singularities the parameter space of those conics is
investigated that have only even order contact with the given quartic.
In particular, its irreducible components are described.
\end{abstract}

\section{Introduction}

The investigations of this article are motivated by problems arising in
the study of twistor spaces over the connected sum of three complex
projective planes. In~\cite{kreussler-kurke} it is shown that, under
suitable conditions, such a twistor space is a small resolution of a
Double Solid branched over a real quartic surface with exactly 13 ordinary
nodes. (A Double Solid is a branched double cover of $\PP ^3$;
see~\cite{clemens} for many aspects of these varieties. The properties
of Double Solids used in this paper are contained in~\cite{kreussler}.)
However, it is not known for which quartics and for which resolutions
twistor spaces occur.  One approach to solving this problem is to study
the family of twistor lines. These are smooth rational curves with
normal bundle $\OO(1)^{\oplus 2}$ in the total space of the twistor
fibration. Those curves will be called ``lines'' in the sequel.  The base
of the twistor fibration must be contained in the set of real points in
the parameter space of all lines in the twistor space. Therefore, it
seems to be promising to investigate the parameter space of all lines
in small resolutions of Double Solids or at least something which is,
hopefully, similar enough to this space.

There are some articles on Double Solids e.g. \cite{clemens} or
\cite{tichomirov} but there is nothing (as far as known to the author)
about lines in Double Solids. In \cite{tichomirov} ``lines'' in Double
Solids are studied but Tikhomirov uses an essentially different notion
of lines. (His lines are mapped to double tangents of the branch locus
-- but cf.  Proposition~\ref{geraden-in-DS}.) Furthermore, he assumes
the branch locus to be smooth.

The following Proposition completes the motivation for the investigations
announced in the abstract:

\begin{satz} \label{geraden-in-DS}
Let $\pi:Z\longrightarrow\PP^3$ be a Double Solid branched over a
quartic $B$ that has at most ordinary nodes as singularities. Let
$C\subset Z$ be a line (i.e., a smooth rational curve with normal
bundle $\OO(1)^{\oplus 2}$) outside the singular set of $Z$. Then the
image $\pi(C)$ of $C$ is a conic that has only even order contact with
$B$.
\end{satz}

By this fact the study of the parameter space of all those ``touching
conics'' is motivated. In this article we will study the parameter space
of touching conics by examining its fibration over $\PPV^3$. (This
fibration is given by assigning to each conic the unique plane in which it
is contained). In Section~\ref{conics-plane} the generic fibre of this
fibration is described. Our approach to the description of the space of
touching conics is the detailed study of the {\em reducible} touching conics
which are complanar pairs of double tangents at the quartic $B$.
Section~\ref{section-dt} is devoted to the examination of the space of
double tangents at $B$. These results are used in
Section~\ref{sec-tochin-conic} to determine the monodromy of the space

\[\YF:=\{(\ell,H)\in \Grass(2,4)\times \PPV^3\,|\, \ell \subset H,
      \ell\;\mbox{is double tangent}\}\]

over $\PPV^3 \setminus \Delta$, where $\Delta$ denotes the closed subset
of planes $H$ such that $B\cap H$ is singular. As a corollary we get the
monodromy of the symmetric product $\left(\YF\symm_{\PPV^3}\YF\right)
\setminus \!\mbox{Diag}$. The latter can be regarded as the parameter
space of reducible touching conics. It serves as a ``frame'' within the
parameter space of all touching conics which, finally, is described using
the knowledge on this ``frame''.

{\bf Acknowledgement:} For their ideas and helpful discussions that
essentially influenced my work on this topic, I am most grateful to Bernd
Kreu{\ss}ler and Prof.~Kurke.

\section{General Preliminaries}\label{gen-prel}

First, for completeness, Proposition~\ref{geraden-in-DS} is to be proven,
now.

{\bf Proof of (\ref{geraden-in-DS}):} Let $\omega_Z$ be the canonical
sheaf of $Z$. It holds $\omega_Z \cong \pi^*\OO_{\PP ^3}(-2)$ (cf.
\cite{bpv} Lemma~I.17.1) and by adjunction formula we have

\hspace*{\mathindent}
\begin{tabular}{@{}r@{\quad}c@{\quad}l@{\extracolsep{3em}}l}
\rule{0mm}{4ex}$\omega _C$ & $\cong$ & $\omega _Z \otimes \det N_{C|Z}$
   & hence\\
\rule{0mm}{4ex}${\cal O}_C(-2)$ & $\cong$ & $\pi ^*{\cal O}_{\PP ^3}
  (-2) \otimes {\cal O}_C(2)$ &\\
\rule{0mm}{4ex}${\cal O}_C(-4)$ & $\cong$ & $\pi ^*{\cal O}_{\PP ^3}
  (-2) \otimes {\cal O}_C$
\end{tabular}

or for divisors on $C$

\[-4[pt] = \pi ^*(-2H)\cdot [C]\qquad \mbox{($H$ -- a hyperplane section
                                                   in $\PP ^3$)} \]

and therefore

\begin{eqnarray*}
\pi _*\left(-4 \left[pt\right]\right) & = & \pi _*\left( \pi ^*
   \left( -2H\right)\cdot \left[C\right]\right) \\
& = &(-2H)\cdot \pi _*\left(\left[C\right]\right) \qquad\qquad\mbox{by
   projection formula and finally}\\
2[pt] &=& H \cdot \pi_*\left( \left[C\right]\right).
\end{eqnarray*}

Thus we get

\[ \pi _*\left(\left[C\right]\right) = 2H^2. \]

Now, if $\pi(C)$ were of degree one in $\PP^3$ and $\pi |_C: C
\rightarrow \pi(C)$ of degree two we would have

\[ N_{C|Z} = N_{\pi^{-1}(\pi(C))|Z} \stackrel{!}{\cong}
             \pi^*N_{\pi(C)|\PP ^3}  \cong
             \pi^*\left( \OO_{\pi(C)}\left( 1\right)^{\oplus 2}\right)
  \cong \OO_C(2)^{\oplus 2} \]

in contradiction to $N_{C|Z} \cong \OO_C(1)^{\oplus 2}$. Therefore
$\pi(C)$ is a rational curve of degree two and $\pi|_C:C\rightarrow
\pi(C)$ is an isomorphism, i.e., $\pi(C)$ is a smooth conic.
Furthermore, $\pi^{-1}(\pi(C))$ must split into two irreducible
components if not $\pi(C)\subset B$.

Now, let $\pi(C)\not\subset B$, $\spec\,A = U \subset \pi(C)$ a suitable
open subset, $f$ the equation of $B$ restricted to $U$, and $A_Z:=
\raisebox{0.25ex}{$A[T]$}\big/ \raisebox{-0.25ex}{$(T^2-f)$}$. Then
$\spec\, A_Z=\pi^{-1}(U) \subset \pi^{-1}(\pi(C))$ and $\spec\, A_Z$ is
reducible if and only if $(T^2-f)$ is reducible in $A[T]$, i.e., if and
only if there is an $g\in A$ with $f=g^2$. This proves the proposition.
\hspace*{\fill} $\Box$

\subsection{Parameter space of conics}\label{psc}

We will construct the parameter space of conics in $\PP ^3$.  The space
of ``touching conics'' will be contained within this parameter space.
Let $\PPV ^3$ be the space of planes in $\PP ^3$ and $\cal S$ the
universal subbundle over $\PPV ^3=\Grass(3,4)$. Then $P:=\PP (Sym^2
{\cal S}^{\vee})$ is the parameter space of all conics in $\PP ^3$.
(Every conic can be given by a plane and a symmetric form of degree two
in this plane. On the other hand, every conic determines a unique plane
which it sits in and in that plane a symmetric 2-form which is unique
up to multiplication by scalars.  Even for a double line there is a
unique plane in which it is contained. It is determined by the
non-reduced subscheme structure of the double line.) The projection
$p:P\longrightarrow \PPV^3$ assigns to each conic the unique plane
which it is contained in.

There is a universal family over $P$, constructed as follows: Let
$H:=\PP(p^*{\cal S})$ be the pull-back of the universal plane over
$\PPV ^3$, $\tau:H\rightarrow P$ the projection, and $\OO_H(1)$ the
relative tautological bundle of $H$ over $P$. Then there is a
distinguished section in $(\tau^* \OO_{P|\PPV ^3}(-1))^{\vee} \otimes
\OO_H(2)$, for it is

\begin{eqnarray*}
\left(\tau ^* \OO_{P|\PPV ^3}\left(-1\right)\right)^{\vee} \otimes
                     \OO_H\left(2\right)
&=&
Hom\left( \tau^* \OO_{P|\PPV ^3}\left(-1\right),
           \left( \OO_H\left( -1\right) ^{\vee} \right)^{\otimes 2} \right) \\
&=&
Hom\left( \tau^*\OO_{P|\PPV ^3}\left(-1\right),
               Sym^2 \left( \OO_H\left( -1\right) \right) ^{\vee} \right).
\end{eqnarray*}

and there are canonical injections of vector bundles over $H$:
\[  {\cal O}_H(-1)\hookrightarrow \tau^*p^*{\cal S}.  \]
\[ \tau^* \OO_{P|\PPV ^3}(-1) \hookrightarrow
   \tau^*p^*(Sym^2 {\cal S}^{\vee}) = \tau^*Sym^2 (p^*{\cal S})^{\vee}.  \]

The distinguished section is given by the composition

\[ \tau^*\OO_{P|\PPV ^3}(-1)\hookrightarrow
   \tau^*Sym^2 p^*{\cal S}^{\vee} \longrightarrow
    Sym^2 {\cal O}_H(-1)^{\vee}. \]

The universal family over $P$ is the zero locus of this section.

The above construction shows that there is a natural projection from
the parameter space of ``touching conics'' onto $\PPV^3$, assigning to
each conic the plane which it sits in. Section~\ref{conics-plane} is
devoted to the study of the fibres of this projection.

\subsection{13-nodal quartics}

As outlined in the introduction, the focus of the present investigations lies
on the study of quartic surfaces with exactly 13 ordinary double points.
Those quartics are extensively studied in \cite{kreussler}. The results
needed in the sequel are to be summarised here.

\begin{satz}\label{quartik13}
Every real\/\myfootnote{with respect to the standard real structure of
$\PP^3_{\!\!\complex}$} quartic surface $B$ with exactly one real point and 13
ordinary double points can be defined by an equation of the form

\begin{eqnarray*}
   F&=&x_3^2f_2+2x_3L_0L_1L_2 +f_2^2-f_2\left(L_0^2+L_1^2+L_2^2\right) +
       L_0^2L_1^2+L_0^2L_2^2+L_1^2L_2^2 \\
    &=& \frac{1}{4}\left( Q^2 - E_1E_2E_3E_4\right) \\
\end{eqnarray*}

where $E_1 = x_3 - L_0 - L_1 - L_2$, $E_2 = x_3 + L_0 + L_1 - L_2$,
$E_3 = x_3 + L_0 - L_1 + L_2$, $E_4 = x_3 - L_0 + L_1 + L_2$, $f_2 =
x_0^2+x_1^2 +x_2^2$, $Q = 2f_2 +x_3^2 - L_0^2 - L_1^2 - L_2^2$ and
$L_j = \sum_{i=0}^2 a_{ij}x_i$ such that $f_2 - L_j^2$ $(j=0,1,2)$ defines
three smooth conics with 12 different intersection points.

If, moreover, the planes $E_i$ are real then the quadratic forms $f_2 -
L_j$ are positive definite. If the forms $L_j$ are mutually linearly
independent then $F$ defines a real quartic with exactly 13 ordinary double
points and exactly one real point $P= (0:0:0:1)$. Each of the six lines
$E_i = E_j = 0$ intersects the quadric $Q$ in two different points which
form a conjugate pair of double points of the quartic.\proved
\end{satz}

{\bf Remark:} Those quartics $B$ satisfying all the conditions of the above
proposition are just the quartics which  generically occur in connection
with twistor spaces as mentioned in the Introduction (cf. \cite{kreussler}
and \cite{kreussler-kurke}).

\begin{lemma}\label{branch-sextic}
The projection of $B$ (as above -- with real planes $E_i$) from $P$ onto
the plane $x_3=0$ defines a double cover of $\PP^2$ which is branched along
the sextic $\ti{S} = (f_2 - L_0^2) (f_2 - L_1^2)(f_2 - L_2^2)$ in $\PP ^2$. The
conic $f_2=0$ touches $\ti{S}$ in six smooth points.\proved
\end{lemma}

\section{Conics touching a plane quartic}\label{conics-plane}

In this section\myfootnote{Many ideas of this section are already contained
in \cite{salmon}.}, for a given quartic curve $B$ in $\PP^2$, the set
of conics that have only even order contact with
$B$ (but are different from double lines) is investigated. Those conics will
be called {\em touching conics} in the sequel.

Throughout this section, let $B\subset \PP ^2$ be an irreducible
quartic curve given by a form  $F$ of degree four which has at most
one ordinary node as its only singularity.

Then the following lemma holds.

\begin{lemma}\label{Strukt-ber-KS}
The set of all touching conics is the union of one-parameter-families.
In each family the elements are mutually different. If $B$ is
smooth, the families  are  pairwise disjoint. If $B$ is
singular, two families can only intersect in a reducible conic
which consists of two lines  both containing the singular point of $B$.
\end{lemma}

\begin{proof} The quadratic form $U$ defines a touching conic if and only
if there exist two further quadratic forms $V$ and $W$ such that

\[  F = UW - V^2. \]

Since

\[ UW - V^2 =U(\lambda ^2U+2 \lambda V +W) - (\lambda U + V)^2 \quad
                    \lambda \in \complex, \]

by $U$, $W$, and $V$ a whole one-parameter-family of touching conics is
given:

\begin{equation} \label{epf}
 \lambda ^2U + 2 \lambda \mu V + \mu ^2 W \quad (\lambda : \mu ) \in
            \PP ^1 .
\end{equation}

All conics in such a family are different from each other.  Otherwise
we would have

\[ F= U(\lambda^2 U +2\lambda V + W) - (\lambda U + V)^2 =
      U(\mu^2 U +2\mu V + W) - (\mu U + V)^2  \]

with $\lambda^2 U +2\lambda V + W  = \mu^2 U +2\mu V + W$ and $\lambda
\ne \mu$. Consequently

\[ (\lambda U + V)^2 = (\mu U+ V)^2 \qquad\mbox{hence}\]
\[ \lambda U+ V = -(\mu U+V) \quad \mbox{thus} \quad
   V= \frac{\lambda + \mu}{2} U  \qquad\mbox{and finally}\]
\[ F=U \left( W-\left( \frac{\lambda +\mu}{2}\right)^2 U\right), \]

so that the quartic $B$ would be reducible in contradiction to the
assumption.

Suppose $U$ is contained in two different families~(\ref{epf}), i.e.,
there exist $V$, $W$, resp. $V'$ and $W'$ satisfying
\[ F= UW -V^2 =UW'-V'^2 \qquad \mbox{and therefore}  \]
\[ U(W-W') =V^2 - V'^2 =(V+V')(V-V'). \]

If $U$ is irreducible it follows $U\,|\,(V+V')$ or $U\,|\,(V-V')$ and
consequently  $V'= \pm (\lambda U-V)$ and $W'=\lambda ^2U-2\lambda V+W$,
i.e., $(U,V,W)$ and $(U,V',W')$ yield the same family. If $U$ is reducible
and neither $U\,|\,(V+V')$ nor $U\,|\,(V-V')$ holds, then the intersection
of the two lines of $U$ must be a point of $B$ which is necessarily
singular. If, finally, $U$ were a double line, $B$ would have an equation
of the form $F =L^2W-V^2$ and therefore would have at least two singular
points if the line $L=0$ intersects $B$ in two points or a cusp in the
intersection point of $B$ with $L=0$. This completes the proof.
\end{proof}

\begin{lemma} \label{epf-verh}
Let $B$ be a quartic with exactly one ordinary node. A conic $U$
consisting of two lines through the node is contained in exactly two
one-parameter-families~(\ref{epf}) that intersect in just this conic.
\end{lemma}

\begin{proof} Let $F= f_1\,x^4 + f_2\,y^4 +\cdots$ be the equation of the
quartic $B$ and let the node of $B$ have the coordinates $(0:0:1)$. In
suitable coordinates $U$ is given by the equation $xy=0$. By assumption,
$F$ is of the form

\begin{equation}\label{coeff-comp}
F=xyW-V^2
\end{equation}

with quadratic forms $V$ and $W$. By comparing the coefficients on both
sides of (\ref{coeff-comp}) one finds exactly two one-parameter-families
$W_\lambda$ containing the conic $U$.
\end{proof}

Now, the number of reducible conics in each family~(\ref{epf}) is to be
determined. Knowing the number of double tangents at a quartic curve,
this permits to count the one-parameter-families of touching conics.

\begin{lemma} \label{epf-zerf-el}
Every one-parameter-family~(\ref{epf}) contains at least five and at most
six reducible elements. The family contains only five reducible conics if
and only if it contains a conic which splits into two lines  intersecting
in a point of $B$ (which is necessarily singular).
\end{lemma}

\begin{proof} A conic is reducible if and only if the determinant of its
matrix of coefficients vanishes. Therefore one has to examine the roots of
the polynomial in $\lambda$: $\det (\lambda ^2 U + 2\lambda V + W)$.  This
polynomial is of sixth degree and, therefore, has at most six roots.  For
proving the lemma it is necessary to investigate the conditions under
which the equation has multiple roots. Since the two triples $(U,V,W)$ and
$(\lambda ^2U + 2\lambda V + W, \lambda U + V, U)$ of conics define the
same one-parameter-family it is sufficient to study under which conditions
$\lambda=0$ is a multiple root. For $\lambda=0$ to be a multiple root,
conditions on the coefficients of $U$, $V$, and $W$ are posed. A simple
calculation shows that a multiple root corresponds to a reducible conic
the lines of which intersect in a point $P\in B$, and that no other conic
of the one-parameter-family then can contain $P$. Therefore, at most one
root of the considered equation can be of higher multiplicity since $B$
was supposed to have at most one singular point. If the equation had a
root of multiplicity greater than two then the quartic $B$, given by the
equation $F=UW-V^2$, would have a cusp.  This proves the lemma.
\end{proof}

The following proposition is proved e.g. in \cite{salmon} or \cite{burau},
the smooth case is also treated in \cite{griffiths-harris} Section~4.4.

\begin{satz}\label{Anz-DT}
Let $B$ be an irreducible quartic curve. If $B$ is smooth then $B$ has 28
double tangents\/\myfootnote{resp. lines that have fourth order contact
with $B$}. If $B$ has one ordinary node as its only singularity then $B$
has 22 double tangents.
\end{satz}

{\bf Remark:} The Pl\"ucker formulas seem to contradict the second part of
the Proposition. One must, however, take into account that lines passing
through the node and touching the quartic in an other point are not counted
by the Pl\"ucker formulas. There are six of those lines as one finds for
instance in \cite{salmon}.

\begin{satz}\label{Anz-epf}
Let $B$ be an irreducible quartic curve. If $B$ is smooth then the set of
touching conics splits into 63 mutually disjoint
one-parameter-families~(\ref{epf}). If $B$ has an ordinary node as its
only singular point then the set of touching conics is divided into 16
families~(\ref{epf}) each of which is disjoint from all other families,
and 15 pairs of families that intersect in exactly one conic which
consists of two lines intersecting in a point of $B$.
\end{satz}

\begin{proof} A smooth quartic curve has 28 double tangents, hence, there
are 378 pairs of double tangents, i.e., 378 reducible touching conics.
Thus, according to Lemma~\ref{Strukt-ber-KS}, there are 63
families~(\ref{epf}) of touching conics and every touching conic is
contained in one of these families.

In the singular case there are:
\begin{enumerate}
\renewcommand{\labelenumi}{\alph{enumi})}
\item 15 pairs of double tangents intersecting in the node of $B$
\item 96 pairs one line of which contains the node
\item 120 pairs no line of which contains the node.
\end{enumerate}

According to Lemma~\ref{epf-verh} each of the pairs "a)" is contained in
two families each of which contains four pairs "c)" (cf.
Lemma~\ref{epf-zerf-el}).  Thus the pairs "a)" and "c)" spread over 30
families which contain exactly these reducible conics. These families
intersect pairwise as stated. The remaining 96 pairs "b)" must be
contained in some families~(\ref{epf}), as well. These families do not
intersect any other family and each contains six of the 96 pairs "b)"
(cf.  Lemma~\ref{Strukt-ber-KS} and Lemma~\ref{epf-zerf-el}). Therefore
the 96 pairs "b)" generate 16 further families.
\end{proof}

\begin{lemma}\label{DT-in-EPF}
Let $B$ be a smooth quartic curve.
\begin{enumerate}
\renewcommand{\labelenumi}{\alph{enumi})}
\item The double tangents occurring in the reducible conics in one
      one-parameter-family~(\ref{epf}) are mutually different.
\item Let $ab$, $cd$, and $eh$ be reducible conics contained in the same
      one-parameter-family~(\ref{epf}) ($a$, $b$, $c$, $d$, $e$, and $h$
      linear forms). Then the double tangent $e$ does not occur in any
      reducible conic of that one-parameter-family in which the conics $ac$
      and $bd$ are contained.
\end{enumerate}
\end{lemma}

\begin{proof} a) If there were reducible conics $ab$ and $ac$ in the same
one-parameter-family then the equation $F$ of $B$ could be written in the
form $F = a^2bc - V^2$ with a quadratic form $V$. The intersection points
of the conic $V$ with the line $a$ then would be singular points of $F$.

b) The quartic $F$ may be written in the form $F = abcd - V^2$ with a
quadratic form $V$. Since $eh$ is contained in the one-parameter-family
spanned by the conics $ab$ and $cd$ there is a $\lambda \ne 0$ such that
$ef = \lambda^2 ab + 2\lambda V +cd$. Now, by writing $F$ in the form $F=
ac\cdot bd-V^2$ one finds the one-parameter-family containing $ac$ and
$bd$ to be $\mu^2 ac + 2\mu\lambda V +\lambda^2bd$. Suppose $eg$ is
contained in the one-parameter-family of $ac$ and $bd$, i.e. there is a
$\mu \ne 0$ such that $eg = \mu^2 ac + 2\mu V +bd$. This yields

\[ V = \frac{1}{2\lambda}\left( ef - \lambda ^2 ab - cd \right) =
       \frac{1}{2\mu}    \left( eg - \mu^2 ac -bd\right)  \]

Hence $e(\mu f - \lambda g) = (\lambda b - \mu c)(\lambda \mu a - d)$, i.e.
the linear form $e$ divides one of the two forms $(\lambda b - \mu c)$ or
$(\lambda \mu a - d)$. Thus the line defined by $e$ contains one of the
intersection points of $b$ and $c$ or of $a$ and $d$. In both cases the
common point of the three lines is a singular point of $B$.
\end{proof}

\section{The parameter space of double tangents}\label{section-dt}

Let $B \subset \PP^3$ be a quartic surface with ordinary double points
as its only singularities which is given by an equation of the form

\begin{equation}\label{quartik-equ}
      g_1 g_3 - g_2^2 = 0
\end{equation}

where $g_i$ are homogeneous of degree $i$. Let $(x_0:\ldots :x_4)$ be
homogeneous coordinates on $\PP^4$ and let $\PP^3 \subset \PP^4$ as the
hyperplane $x_4 = 0$. Then

\begin{equation}\label{kubik}
      x_4^2g_1 + 2x_4 g_2 +g_3 = 0
\end{equation}

defines a cubic $K$ in $\PP^4$ with $P:= (0:0:0:0:1)\in K$. Consider the
projection from $P$ onto the hyperplane $x_4 = 0$. The extension
$\ti{K}\stackrel{\pi}{\longrightarrow} \PP^3$ of this map to the blow-up
$\ti{K}$ of $K$ in $P$ is a partial small resolution
of the double solid $Z_0$ branched over $B$ (cf. \cite{kreussler-pre}). The
$\pi$--fibre of any point $x \in \PP^3$ with $g_1(x)=g_2(x)=g_3(x)=0$ (i.e.
the singular points of $B$ in the plane $g_1=0$) is
just the strict transform of the line through $P$ and $x$, all other fibres
consist of at most two points. In particular, there are exactly six lines
through $P$ in $K$ (if the three surfaces $g_i = 0$ intersect properly)
namely the six lines that are contracted by $\pi$. The only singularities
of $K$ are the preimages of double points of $B$ not contained in the
plane $g_1 = 0$. These are ordinary double points, as well.

\begin{lemma}
$\pi: K\setminus \{P\} \longrightarrow \PP^3$ maps lines in $K\setminus
\{P\}$ onto lines that have even intersection with $B$\myfootnote{We will
call those lines simply double tangents, i.e. lines with fourth order
contact will be called double tangents, too.} or which are contained in $B$.
If such a line is contained in $B$ the it passes through a singular point
of $B$ which is contained in the plane $g_1=0$. Moreover, such a line
has one point of higher order intersection with the cubic $g_3=0$ or the
line is contained in this cubic.
\end{lemma}

\begin{proof} Let

\begin{morph}
   L:&\PP^1 & \lhookrightarrow & \PP^4 \\
   &(s:t)   & \longmapsto      & (l_0(s,t):\ldots:l_4(s,t))
\end{morph}

be the parameter representation of a line $L$ in $\PP^4$, which is
contained in $K$. Thus the equation~(\ref{kubik}) restricted to $L$
vanishes, i.e.

\[ \left.g_3\right|_L = -l_4^2 \left.g_1\right|_L -2l_4 \left.g_2\right|_L
   \qquad\mbox{thus} \]

\begin{eqnarray*}
   \left. g_1\right|_L\cdot\left. g_3\right|_L - \left. g_2\right|_L^2 &=&
     -l_4^2 \left. g_1\right|_L^2 -2l_4 \left. g_1\right|_L\cdot \left.
     g_2\right|_L - \left. g_2\right|_L^2 \\
   &=& - \left( l_4 \left. g_1\right|_L + \left. g_2\right|_L
      \right)^2
\end{eqnarray*}

Hence, equation~(\ref{quartik-equ}) restricted to $\pi(L)$ is a complete
square and therefore the image of $L$ in $\PP^3$ is a line with even
intersection with $B$ or is contained in $B$. If $\pi(L)$ is contained in
$B$ then $L$ is contained in the ramification locus of $\pi$. Hence, the
plane spanned by $L$ and $P$ intersects $K$ in a cubic curve which consists
of $L$ counting twice and a line through $P$. The lines through $P$ in $K$
are the lines connecting $P$ with the singular points of $B$ in the plane
$g_1=0$. Moreover, if $\pi(L)\subset B$ then $(-\l_4\cdot g_1|_L + g_2|_L)$
must vanish and, hence, $g_3|_L = \l_4^2\cdot g_1|_L$.
\end{proof}

\hspace*{\fill}
\begin{minipage}{0.95\textwidth}\sloppy
{\bf Remark:} In the case that we are particularly interested in -- $B$ is
given by an equation as described in Proposition~\ref{quartik13} and $g_3$
is the product of three of the linear forms $E_i$, say $g_3=E_2E_3E_4$.
$Q$ and $E_i$ cannot have common real zeros for there is only one real
point on $B$ which is outside the quadric $Q=0$. Therefore the planes
$E_i=0$ intersect the quadric $g_2:=Q=0$ along smooth conics.
Consequently, there is no line which is contained in $B$ and in the cubic
$g_3=0$.\par\smallskip

A line $\ell$ (not contained in $g_3=0$) that has a point $P_a$ of higher
order intersection with the cubic $g_3=0$ must meet the intersection of
two of the three planes the cubic consists of. On the other hand, no
common point of three of the four planes $E_i$ is a point of $B$ by
Proposition~\ref{quartik13}. In particular, if $\ell\subset B$ then $P_a$
is not contained in the plane $g_1=0$. Hence, if $\ell$ is the image of a
line in $K$ which is contained in $B$ then $\ell$ must pass through two
singular points of $B$, namely $P_a$ and the the singular point $P_b$ of
$B$ in the plane $g_1$ which $\ell$ must pass through by the above lemma.
\par\smallskip

Therefore, in our case only those lines in $B$ (if any) can appear as
images of lines in $K$ that pass through two double points of $B$. But
lines through two double points are double tangents unless they are
contained in $B$ and hence, lines in $B$ through two double points
necessarily appear in the parameter space of double tangents.
\end{minipage}

\begin{lemma}
Every double tangent of $B$ that is not contained in the plane $g_1 =0$ is
the image under $\pi$ of a line in $K$.
\end{lemma}

\begin{proof}  Let

\begin{morph}
   L':&\PP^1 & \lhookrightarrow & \PP^3 \\
      &(s:t)   & \longmapsto      & (l_0(s,t):\ldots:l_3(s,t))
\end{morph}

be the parameter representation of a line in $\PP^3$. If $L'$ is contained
in the surface $g_3=0$ then there is nothing to show; so, let $L'$
not be contained in $g_3 = 0$. If $L'$ is a double tangent then there
exists a quadratic form $q$ on $L'$ such that

\[ \left. g_1\right|_{L'}\cdot\left. g_3\right|_{L'} - \left. g_2
   \right|_{L'}^2 = - q^2 \qquad\mbox{and therefore} \]
\[ \left. g_1\right|_{L'}\cdot\left. g_3\right|_{L'} =
   \left( \left. g_2\right|_{L'} + q\right)\left(  \left. g_2\right|_{L'}
    - q\right) \]

Hence there is a linear form $l_4(s,t)$ on $L'$ such that $l_4 \left.
g_1\right|_{L'} = - \left( \left. g_2\right|_{L'} \pm q\right)$. By
eventually replacing $q$ with $-q$ we can assume that $l_4 \left.
g_1\right|_{L'} = - \left( \left. g_2\right|_{L'} + q\right)$. Then

\begin{eqnarray*}
   \lefteqn{\left. g_1\right|_{L'}\cdot \left( l_4^2 \left. g_1\right|_{L'} +
      2l_4 \left. g_2\right|_{L'} +  \left. g_3\right|_{L'}\right) =}
      \hspace{7em}\\[1ex]
   &=& \left( \left. g_2\right|_{L'} + q\right)^2 +2l_4\left.
       g_1\right|_{L'}\cdot\left. g_2\right|_{L'} + \left( \left.
       g_2\right|_{L'} + q\right)\left(  \left. g_2\right|_{L'} - q\right)\\
   &=& \left( \left. g_2\right|_{L'} + q\right)\left( \left.
      g_2\right|_{L'} + q +  \left. g_2\right|_{L'} - q\right) + 2l_4
      \left. g_1\right|_{L'}\cdot\left. g_2\right|_{L'}\\
   &=& 2 \left. g_2\right|_{L'}\left(\left. g_2\right|_{L'} + q + l_4\left.
          g_1\right|_{L'}\right)\\
   &=& 0
\end{eqnarray*}

Therefore, since $\left. g_1\right|_{L'}\not\equiv 0$, $l_4^2 \left.
g_1\right|_{L'} + 2l_4 \left. g_2\right|_{L'} +  \left. g_3\right|_{L'}=0$
and hence $(l_1(s,t):\ldots :l_4(s,t))$ ($(s:t) \in \PP^1$)
defines a line in $\PP^4$ which is contained in $K$ and the image under
$\pi$ of which is $L'$.
\end{proof}

\begin{lemma}\label{dt-preimage}
Let $L$ be a double tangent of $B$ which is not contained in $g_1=0$. If
$L'$ and $L''$ are different lines in $K$, which are mapped onto $L$ then
$L$ contains one of the points $x$ with $g_1(x)=g_2(x)=g_3(x)=0$. In this
case $L'$ and $L''$ are the only lines in $K$ that are mapped to $L$.
\end{lemma}

\begin{proof} Consider the plane $H$ in $\PP^4$ spanned by $L$ and $P$.
Then $H$ is not contained in $K$. If $H$ were contained in $K$ then $L$
would be contained in the locus $g_3 = 0$ in $\PP^3$. $H\subset K$ then
implies $g_1=g_2=0$ on $L$. But $L$ was supposed to be not contained in
$g_1=0$.  Therefore $K\cap H$ is a plane cubic curve that contains $L'$
and $L''$.  Hence, this plane cubic must split into three components all
of which are lines. One of these lines must contain $P$, but the lines
through $P$ are just the lines $\overline{Px}$ with $x \in \PP^3$ and
$g_1(x)=g_2(x)=g_3(x)=0$. On the other hand $x$ is in $L$ which proves the
lemma.
\end{proof}

The projection $\pi$ induces a morphism $\bar{\pi}$ between Grassmannians of
lines of $\PP^4$ and $\PP^3$:

\begin{morph}
      \bar{\pi} :& \Grass(2,5)\setminus \{\mbox{lines through $P$}\}&
         \longrightarrow & \Grass(2,4) \\
      &\cup &&\cup \\
      &\Fano(K)\setminus\{\mbox{lines through $P$}\} & \longrightarrow &
         \left\{
         \parbox{0.35\textwidth}{lines with even intersection with $B$ and
	    lines contained in $B$}
         \right\}
\end{morph}

where $\Fano(K)$ denotes the set of lines in $K$. The above lemmata suggest
that $\bar{\pi}$ defines a birational map onto the set of bitangents of $B$
in $\Grass(2,4)$. This will be shown later.

\subsection{Lines in a nodal cubic threefold and bisecants of a space curve}

Let $K \subset \PP^4$ be a cubic hypersurface with only ordinary double
points as singularities. Let \POne{}  be one of the double points, $H$ a
hyperplane not containing \POne{}, and let $Q'$ be the intersection of the
tangent cone of \POne{}  with the hyperplane $H$. Finally, let $S = K \cap
Q'$. Assume the curve $S$ to have only ordinary double points as
singularities. (Under the above assumptions on $K$ this is, in fact,
always true, cf. \cite{finkelnberg}.) The singularities of $S$ are just
the images of the double points of $K\setminus \{\POne{} \}$. Denote by
$p$ the projection $p:K\setminus \{\POne{} \} \longrightarrow H$. Then the
following proposition holds:

\begin{satz}\label{bisecs-lines}
Let $\ell$ be a line in $K$. If $\POne{} \in \ell$ then $\ell$ is a line
$\bar{\POne{} x}$ with $x \in S$ and, conversely, every $x \in S$ defines
a line in $K$ through \POne{}. A line $\ell$ not containing \POne{} is
mapped onto a line $\bar{\ell} = p(\ell) \subset H$ not contained in $Q'$
and either connecting two points of $S$ or being a tangent of $S$ at a
smooth point of $S$ or a tangent of $Q'$ at a non-smooth point of $S$.
Conversely every such line $\bar{\ell}$ is the image of a line in $K$.
\end{satz}

\begin{proof}\newcommand{\lbar}{\mbox{$\bar{\ell}$ }}
Here, a sketch of the proof of Finkelnberg (cf.  \cite{finkelnberg}) is to
be given. The case where $\POne \in \ell$ is obvious. Let \lbar be an
arbitrary line in $H$ and let $V$ be the plane in $\PP^4$ spanned by
$\lbar$ and \POne{}. Noting that $Q' \cong \PP^1 \times \PP^1 \subset
\PP^3$ and $S \in \left|\OO_{Q'}(3,3)\right|$ the following cases occur:

\begin{itemize}
\item {\it $\lbar \cap S = \emptyset$ and \lbar intersects $Q'$ transversally.}
      Then $V \cap K$ is an irreducible plane cubic with an ordinary double
      point in \POne{}.
\item {\it $\lbar \cap S = \emptyset$ and \lbar is tangent to $Q'$.} Then $V
      \cap K$ is an irreducible plane cubic with a cusp in \POne{}.
\item {\it $\lbar \cap S = \{T\}$ and \lbar intersects $Q'$ transversally.}
      Then $V \cap K$ consists of the line $\bar{\POne{} T}$ and a conic
      intersecting $\bar{\POne{} T}$ transversally. The conic must pass through
      \POne{}  and so can not split into two lines since one of the lines were
      a line through \POne{}  and would intersect $S$ in a point different
      from $T$. If $\bar{\POne{} T}$ were tangent to the conic then $V \cap K$
      had a cusp which is impossible since \lbar intersects $Q'$ in two
      points.
\item {\it $\lbar \cap S = \{T\}$, $T$ is a non-singular point of $S$ and
      \lbar is tangent to $Q'$ but not to $S$.} Then  $V \cap K$ consists of
      a smooth conic through \POne{}  and the line $\bar{\POne{} T}$ which is
      tangent to the conic (by a similar argument as above).
\item {\it $\lbar \cap S = \{T_1,T_2\}$ with $T_1 \ne T_2$ and \lbar is not
      contained in $Q'$.} Then
      $V \cap K$ contains the two lines $\bar{\POne{} T_1}$ and
      $\bar{\POne{} T_2}$. Even if one or both points $T_i$ are singular
      points of $S$ none of these lines counts twice since this is only
      possible if \lbar is tangent to $Q'$. Therefore $V \cap K$ contains a
      third line not through \POne{}.
\item {\it \lbar is tangent to $S$ in a smooth point $T$ of $S$.} Then
      $V \cap K$ splits into the line $\bar{\POne{} T}$ with multiplicity two
      and a second line that does not pass through any singular point of
      $K$. (To see that the first line carries multiplicity two move \lbar
      a bit.)
\item {\it \lbar is tangent to but not contained in $Q'$ and meets $S$ in a
      singular point $T$.} Then
      $V \cap K$ consists of the line $\bar{\POne{} T}$ with multiplicity two
      (move \lbar a bit!) and a second line. The second line can not pass
      through \POne{} since then it would intersect $S$ in a point
      different from $T$.
\item {\it \lbar is contained in $Q'$.} Then \lbar meets $S$ in three (not
      necessaryly different) points. Each of these points is the intersection
      point of a line  through \POne{} in $K$ with $S$, multiple points
      corresponding to multiple lines.
\end{itemize}
Therefore, only those lines mentioned in the proposition can be images of
lines in $K$ not through \POne{} and each of them determines such a line in
$K$. This proves the proposition.
\end{proof}

\subsection{The Fano scheme of lines on the cubic threefold}

{}From now on it becomes convenient to make use of the special type of the
quartic $B$ as described in Proposition~\ref{quartik13}: Let $B$ be a real
quartic with exactly 13 ordinary double points such that (using the notation
of Proposition~\ref{quartik13}) the linear forms $E_1,\ldots ,
E_4$ are real (i.e. have real coefficients). Take $g_1 = E_1$,
$g_3 = E_2E_3E_4$ and $g_2 = Q$ to achieve the form $g_1g_3 - g_2^2$, i.e.
$B$ is given by an equation of the form

\begin{eqnarray*}
   F&=&x_3^2f_2+2x_3L_0L_1L_2 +f_2^2-f_2\left(L_0^2+L_1^2+L_2^2\right) +
       L_0^2L_1^2+L_0^2L_2^2+L_1^2L_2^2 \\
    &=& \frac{1}{4}\left( Q^2 - E_1E_2E_3E_4\right) \\
    &=& \frac{1}{4}\left( g_2^2 -g_1g_3 \right) \\
\end{eqnarray*}

With the notation of the previous paragraph an letting $H\subset \PP^4$ be
the plane $x_4=0$, the following lemma holds.

\begin{lemma}
The curve $S = K \cap Q'$ consists of three components, each a smooth conic
in one of the planes $E_i = 0$ ($i=2,3,4$). Each two components intersect in
two points.
\end{lemma}

\begin{proof} First, the equation of $Q'$ is determined: Let $T$ be a
point in $H$. Then $T \in Q'$ if and only if the line $\bar{\POne{} T}$ is
contained in the tangent cone of $K$ in \POne{}  which is the case if and
only if this line has third order contact with $K$ in \POne{}.  A simple
calculation then shows that $T\in Q'$ if and only if $T$ is contained in
the zero locus of the quadratic form $\left(x_3 + L_0 +L_1 +L_2\right)^2 -
4f_2$ (notation as in Proposition~\ref{quartik13}). Since $S= K \cap Q' =
(K\cap H) \cap Q'$ and since $K \cap H$ is the union of the three planes
$E_i = 0$ ($i=2,3,4$), the components of $S$ are the three conics $Q'\cap
\{E_i=0\}$ (see also last remark in \cite{kreussler-pre}). These three
conics are given by the equations

\begin{eqnarray}
f_2 - L_2^2 & = & 0 \qquad\mbox{in the plane $E_2=0$} \nonumber\\
f_2 - L_1^2 & = & 0 \qquad\mbox{in the plane $E_3=0$} \label{S-equ}\\
f_2 - L_0^2 & = & 0 \qquad\mbox{in the plane $E_4=0$} \nonumber
\end{eqnarray}

In the plane $x_3 = 0$, these equations define smooth conics  (by
Proposition~\ref{quartik13}) and since none of the planes $E_i=0$ contains
the point $(0:0:0:1)$ they define smooth conics in the planes $E_i=0$, as
well. Another short calculation shows that the quadrics $Q$ and $Q'$
coincide on the lines $E_i = E_j = 0$. Therefore, by the assumptions on
the quartic and Proposition~\ref{quartik13}, these lines intersect $Q'$ in
two different points.
\end{proof}

Denote by $\Bisec(S)\subset \Grass(2,4)$ the closure of the set of
bisecants of $S$. By the above lemma, $S$ splits into 3 irreducible
components which will be denoted by $S_i,\quad i=1,2,3$. Let $B_{ij}
\subset \Bisec(S) \quad (1\le i\le j \le 3)$ be the closure of the set of
bisecants that connect $S_i$ with $S_j$. $B_{ii}$ then consists of the
bisecants and the tangents of $S_i$ whereas $B_{ij}$ contains the lines
connecting different points of $S_i$ and $S_j$ together with all tangents
at $Q'$ in the to intersection points of $S_i$ and $S_j$. The $B_{ij}$ are
the six irreducible components of $\Bisec(S)$, for they are irreducible,
cover all of $\Bisec(S)$ and contain open subsets that are disjoint from
all other $B_{ij}$.

Now, using Proposition~\ref{bisecs-lines}, morphisms from $B_{ij}$ to
$\Fano(K)$ are to be constructed which will turn out to induce  a
birational map between $\Bisec(S)$ and $\Fano(K)$.  Let $\,{\cal
U}_{B_{ij}} \subset \PP^3 \times B_{ij}$ be the ''universal line'' over
$B_{ij}$: ${\cal U}_{B_{ij}} = \{(x,\ell)\,|\,x \in \ell \}$. Let
$C_{B_{ij}}$ be the cone from $\{\POne{}\} \times B_{ij} \subset \PP^4
\times B_{ij}$ over $\,{\cal U}_{B_{ij}}$ (where $\PP^3 \subset \PP^4$ as
the hyperplane $x_4=0$):

\[
\newlength{\temp}
\settowidth{\temp}{$x$ is contained in the plane spanned by $\ell$ and \POne{}
                   in $\PP^4 \times \{\ell\}$}
   C_{B_{ij}}:= \left\{ (x,\ell) \in \PP^4 \times B_{ij} \left|\;
                \parbox{0.6\temp}{$x$ is contained in the plane
                spanned by $\ell$ and \POne{}  in $\PP^4 \times \{\ell\}$}
                \right.\right\}
\]

Finally let $K_{B_{ij}}:= C_{B_{ij}}\cap \left(K\times B_{ij}\right)$. Every
fibre of $K_{B_{ij}}$ over $B_{ij}$ splits into lines at least two of which
contain \POne{}  (possibly one line counted twice if $\ell \in B_{ij}$ is a
tangent at $S$).

\[\renewcommand{\arraystretch}{1.5}
\begin{array}{ccccl}
      \PP^3 \times B_{ij} & \subset & \PP^4 \times B_{ij} & \longrightarrow
         & B_{ij}\\
      \cup & & \cup\\
      {\cal U}_{B_{ij}} & \subset & C_{B_{ij}} & \supset &
         \{\POne\}\times B_{ij} = \mbox{vertex of the cone}\\
      & & \cup\\
      & & K_{B_{ij}} & \supset & G_{B_{ij}}\\
      & & \cup\\
      \left(\left(S_i \cup S_j \right)\times B_{ij} \right) \cap
         {\cal U}_{B_{ij}} & \subset & M_{B_{ij}} & = & \mbox{cone over}
         \quad \left(\left(S_i \cup S_j \right)\times B_{ij} \right)\cap
         {\cal U}_{B_{ij}}
\end{array}
\]

Consider now $(S_i \cup S_j)\times B_{ij}\cap {\cal U}_{B_{ij}}$ the fibre
over $\ell \in B_{ij}$ of which consists of the intersection points of
$\ell$ with $S_i \cup S_j$. Let $M_{B_{ij}} \subset \PP^4 \times B_{ij}$
be the cone from $\{\POne{}\} \times B_{ij}$ over $(S_i \cup S_j)\cap
{\cal U}_{B_{ij}}$ which by construction is contained in $K_{B_{ij}}$.
Since a line $\ell \in B_{ij}$ can meet $S$ in three different points if
and only if it meets all three components of $S$ the fibre over $\ell$
consists of one or two lines through \POne{} - one of them possibly
counting twice. Let $G_{B_{ij}}$ be the closure in $\PP^4 \times B_{ij}$
of $K_{B_{ij}} \!\setminus M_{B_{ij}}$.

By Proposition~\ref{bisecs-lines} $G_{B_{ij}}$ is a family of lines in
$\PP^4$ all lying in the cubic $K$. Let $\ell \in B_{ij}$ be an arbitrary
line and let $V$ be the plane in $\PP^4$ spanned by $\ell$ and \POne{}. By
subtracting $M_{B_{ij}}$ from $K_{B_{ij}}$ just those lines in $K\cap V$
are removed that pass through \POne{} and through the intersection points
of $\ell$ with $S_i$ and $S_j$. For each line $\ell \in B_{ij}$ such that
$V\cap K$ contains a line not through \POne{} just this line remains in
$G_{B_{ij}}$.  These $\ell$ are just the lines in $\Bisec(S)$ that do not
meet all three components of $S$. But also for the lines $\ell$ that are
contained in $Q'$ the fibre of $G_{B_{ij}}$ over $\ell$ consists of one
single line.  $\ell$ must meet all three components of $S$ and, therefore,
must be contained in each $B_{ij}$ with $i\ne j$. Subtracting $M_{B_{ij}}$
from $K_{B_{ij}}$ removes just the lines in $V\cap K$ that intersect $S_i$
or $S_j$ leaving the the third line which intersects the third component
of $S$ in $G_{B_{ij}}$. Since $G_{B_{ij}}$ is the {\em closure} of
$K_{B_{ij}} \!\setminus M_{B_{ij}}$ this construction also yields just one
line in the case that $\ell$ meets a singular point of $S$ which is always
the intersection of two components of $S$.

Since $G_{B_{ij}} \longrightarrow B_{ij}$ is a family of lines in $K$
there is a uniquely determined morphism $B_{ij} \longrightarrow
\Fano(K)$.  The images of all $B_{ij}$ cover $\Fano(K)$. For lines in $K$
that do not contain \POne{} this is clear from
Proposition~\ref{bisecs-lines}. But also the lines through \POne{} have
their representation by an element of one of the $B_{ij}$:  Let $T$ be the
intersection of such a line with $S$. There are two lines of the rulings
of $Q'$ through $T$. Since $S_i \in \left|\OO_{Q'}(1,1)\right|$ for all
$i$ non of them can be tangent to $S$. Let $\ell$ be one of these two
lines. $\ell$ intersects all three components of $S$ and is therefore
contained in each $B_{ij}$ with $i \ne j$. Therefore, we can find $i$ and
$j$ such that $T$ is neither contained in $S_i$ nor in $S_j$ and $\ell \in
B_{ij}$. (If $T\in S_m \cap S_n$ then e.g. $i=m$ and $j \ne n$ is a
fitting choice.) As is clear from the above discussion, the morphism
$B_{ij} \longrightarrow \Fano(K)$ then maps $\ell \in B_{ij}$ to the line
$\bar{\POne{}T}\in \Fano(K)$.

Now, each $B_{ij}$ contains an open subset (say $U_{ij}$) on which the
above morphisms to $\Fano(K)$ are injective and such that the images of
different $U_{ij}$ do not intersect in $\Fano(K)$. (The sets of lines that
intersect $S$ in exactly two different nonsingular points will do.)
Therefore the morphisms $G_{B_{ij}} \longrightarrow B_{ij}$ induce a
birational map $\Bisec(S)\ratarrow \Fano(K)$.

\subsection{The components of the space of double tangents}

Denote by $\YO \subset \Grass(2,4)$ the closed subscheme of double
tangents of the quartic $B$. Earlier in this section we constructed a
morphism $\Fano(K)\setminus \{\mbox{lines through $P$}\} \longrightarrow
\YO$. By Lemma~\ref{dt-preimage} this morphism is injective outside the
closed subset of lines in $K$ that meet one of the six lines through $P$.
As this closed subset is not an irreducible component of $\Fano(K)$ (which
is clear from the map $\Bisec(S) \ratarrow \Fano(K)$) its complement is
open und dense. On the other hand, the image of $\Fano(K)\setminus
\{\mbox{lines through $P$}\}$ in $\YO$ contains the set of all double
tangents outside the plane $g_1=0$.

But the double tangents contained in this plane form an irreducible
component of their own: Every line in that plane is a double tangent. Thus
the set of these double tangents is closed in $\YO$ and of dimension two.
On the other hand the dimension of $\YO$ is also two: Consider the
variety

\[  \YF{}:=\{ (\ell,H) \in \Grass(2,4)\times \PPV^3\,|\, \ell\subset H \;
      \mbox{and}\; \ell\in \YO \}    \]

which is a subvariety of the flag variety $F(1,2):=\{ (\ell,H) \in
\Grass(2,4) \times \PPV^3 \,|\, \ell\subset H\}$.  This variety is fibred
over $\PPV^3$ and the fibre over a general $H\in\PPV^3$ is zero
dimensional by Proposition~\ref{Anz-DT}. Hence the dimension of \YF{} is
at least three and since every line is contained in a pencil of planes the
dimension of \YO{} is at least two. The dimension cannot be greater than
two for there is a morphism $\Fano(K)\setminus \{\mbox{lines through
$P$}\} \longrightarrow \YO$ (which is surjective onto the set of double
tangents outside the plane $g_1=0$) and the dimension of $\Fano(K)$ is two
since it is birationally equivalent to $\Bisec(S)$. Therefore \YO{} is the
union of two-dimensional varieties and hence is of dimension two. So, the
closed subvariety of lines in the plane $g_1=0$ has the same dimension as
\YO{} and thus is an irreducible component.

This way we have found a rational map $\Bisec(S)\ratarrow\YO$ which is
birational onto those components of \YO{} which are different from the
component of lines in the plane $g_1=0$. \YO{} must therefore have seven
irreducible components. Four of them consist of all lines in a plane.
These are the sets of lines in the planes $E_i = 0$ $(i=1,\ldots,4)$.  As
the lines in the planes $E_i = 0$ $i=2,3,4$ are contained in $K$ as well
as in the hyperplane $x_4 = 0$ they keep fixed under the map
$\Bisec(S)\ratarrow\YO$ and therefore correspond to the lines in the
$B_{ii} \subset \Bisec(S)$.

Thus, the irreducible components of \YO{} are determined:

\begin{satz}\label{DT-compon}
\YO{} consists of seven irreducible components, namely the four components
with all lines in the planes $E_i$ $(i=1,2,3,4)$ and three further
components corresponding to the components $B_{ij} \subset \Bisec(S)$ with
$i\ne j$.\hspace*{\fill}$\Box$
\end{satz}

The following proposition will be needed in the next section.

\begin{satz}\label{flach}
Let $\Delta' \subset \PPV^3$ be the closed set of planes $H$ such that
$B\cap H$ is a plane quartic which is not smooth and has more or other
singularities than just one ordinary double point. Then $\varphi
:\YF|_{\PPV^3 \setminus \Delta'} \longrightarrow \PPV^3 \setminus
\Delta'$ is flat. The ramification locus of $\varphi$ is just the set of
those $(\ell,H)$ such that $B\cap H$ is a quartic with a node and $\ell$
is a double tangent containing the node. The ramification index in those
points is two. Outside the ramification locus $\varphi$ is a smooth
28-fold cover.
\end{satz}

\begin{proof} \YF{} is contained in the Flag variety $F(1,2)\subset
\Grass(2,4) \times \PPV^3$ which is a $\PP^2$-bundle over $\PPV^3$. For
the proof of the flatness we will determine the Hilbert polynomial of the
fibres of $\YF|_{\PPV^3 \setminus \Delta'}$ over $\PPV^3 \setminus
\Delta'$ considered as subsets of $\PPV^2$.

Let $U\subset \PPV^3$ be a standard open subset (say $U=\{(y_0:\ldots:y_3)
\in \PPV^3\,|\, y_3 = -1\}$) such that $F(1,2)|_U\cong U\times \PPV^2$. In
$\PPV^2$ choose a standard open set $U' = \{(l_0:l_1:l_2) \in \PPV^2\,|\,
l_0 = -1\}$. Then, by substituting $x_3 = x_0y_0+x_1y_2 +x_2y_2$ and $x_0=
x_1l_1 + x_2l_2$ in the equation of $B$ we get a map from $U\times U'$ to
$\PP^4$ -- the set of binary quartic forms -- which associates to each
pair $(\ell,H)\in U\times U'$ the equation of $B$ restricted to $\ell$.
If $C\subset \PP^4$ denotes the closed subset parametrising those quartic
forms that are complete squares then $\YF|_{U\times U'}$ is the preimage
of $C$ under the map $U\times U' \longrightarrow \PP^4$. Therefore, we can
get equations for $Y|_{U\times U'}$ (with its reduced scheme structure) by
pulling back equations defining $C$.

Let $(a:b:c:d:e)\in \PP^4$ correspond to the quartic form $ax^4 + bx^3y +
cx^2y^2 + dxy^3 + ey^4$. Obviously we get the equations for $C$ with the
reduced scheme structure by eliminating $\xi$, $\eta$, and $\zeta$ from
the equations

\begin{equation}\label{triv-equ}
\renewcommand{\arraystretch}{1.5}
\begin{array}{rcl}
  a & = & \xi^2\\
  b & = & 2\xi\eta\\
  c & = & 2\xi\zeta-\eta^2\\
  d & = & 2\eta\zeta\\
  e & = & \zeta^2.
\end{array}
\end{equation}

This is done by the technique of Gr\"obner bases (cf. \cite{cox} for
details): One has to compute a Gr\"obner basis for the above equations
using the lexicographic monomial order induced by the ordering of
variables $(\xi>\eta>\zeta>a>b>c>d>e)$.  The equations defining $C\in
\PP^4$, then, are those which do not contain $\xi$, $\eta$ or $\zeta$
(cf.~\cite{cox}). The following is the Gr\"obner basis computed by
MAPLE~V (the polynomials not containing $\xi$, $\eta$ or $\zeta$ listed
last):

\begin{equation}\label{GB}
\renewcommand{\arraystretch}{1.5}
\begin{array}{ccc}
   \xi ^{2}-a,& 2\,\xi \eta -b,& 2\,\xi \zeta -c+\eta ^{2},\\
   \xi b-2\,a\eta ,& -\eta b-4\,\zeta a+2\,\xi c,& \xi d-b\zeta ,\\
   4\,\xi e+\eta d -2\,\zeta c,& b\zeta -\eta c+\eta ^{3},& 4\,a\eta
      ^{2}-b^{2},\\
   b\eta ^{2}-cb+2\,ad,& 2\,c\eta ^{2}-2\,c^{2}+ bd+8\,ea,& d\eta
      ^{2}-dc+2\,eb,\\
   2\,\eta \zeta -d,& 4\,a\eta c-4\,ab\zeta -b^{2}\eta ,& 2\,a\eta
      d-b^{2}\zeta ,\\
   4 \,a\zeta d-2\,cb\zeta +\eta bd,& dc\eta -2\,\zeta c^{2}+\zeta
      bd+8\,ea\zeta ,& \eta d^{2}-2\,\zeta dc+4\,eb\zeta ,\\
   2\,\eta e-\zeta d,& \zeta ^{2}-e\\[3ex]
   8\,a^{2}d-4\,cba+b^{3},& cb^{2}+2\,abd-4\,ac^{2}+16\,a^{2}e,&
      b^{2}d+8\,bea-4\,acd,\\
   ad^{2}-eb^{2},& bd^{2}-4\,ceb+8\,ead,& d^{2}c-4\,ec^{2}
      +2\,ebd+16\,e^{2}a,\\
   d^{3}-4\,edc+8\,e^{2}b
\end{array}
\end{equation}

One can easily verify that these polynomials form a Gr\"obner basis with
respect to the above monomial order (e.g. by forming S-polynomials and
reducing them with respect to the set of polynomials~(\ref{GB})). Even
easier is it to check that the equations~(\ref{triv-equ}) and the
polynomials~(\ref{GB}) define the same ideal. (The
equations~(\ref{triv-equ}) are among the polynomials~(\ref{GB}) and the
remaining polynomials are easily reduced to zero using the
equations~(\ref{triv-equ}).)

Now, let $H\in U\setminus\Delta'$ be a plane in $\PP^3$. For any
$(\ell,H)$ in the fibre $\varphi^{-1}(H)$ of \YF{} over $H$ we will compute
the length of the local ring $\OO_{\varphi^{-1}(H),(\ell,H)}$. Let

\begin{eqnarray*}
f & := & a_{44}\,x_0^4 +a_{43}\,x_0^3x_1 +a_{42}\,x_0^2x_1^2
   +a_{41}\,x_0x_1^3 +a_{40}\,x_1^4\\
&& +x_2\,(a_{33}\,x_0^3+a_{32}\,x_0 ^2x_1 +a_{31}\,x_0x_1^2 +a_{30}\,x_1^3) \\
&& +x_2^2(a_{22}\,x_0^2 +a_{21}\,x_0x_1 +a_{20}\,x_1^2)\\
&& +x_2^3(a_{11}\,x_0 +a_{10}\,x_1) + a_{0}\,x_2^4
\end{eqnarray*}

be the polynomial defining $B\cap H \subset H \cong \PP^2$. Suppose that $\ell
\subset \PP^2$ is the line $x_0=0$ contained in $U\times U' \subset
F(1,2)$. The equations of \YF{} near $(\ell,H)$ are obtained by
substituting $x_0=x_1l_1+x_2l_2$ in $f$ and substituting the resulting
coefficients in the seven last polynomials of (\ref{GB}). (The local ring
$\OO_{\varphi^{-1}(H),(\ell,H)}$ then is obtained as the factor of
$\complex[l_1,l_2]_{(0,0)}$ by the ideal generated by the equations of
\YF{} we got in this way.)

If $\ell\subset H$ is a double tangent at $B \cap H$ touching $B\cap H$ in
the points $(0:1:0)$ and $(0:0:1)$ then the coefficients in $f$ have to
satisfy: $a_0=0$, $a_{10}=0$, (since $x=0$ is a tangent at $(0:0:1)$),
$a_{40}=0$, $a_{30}=0$ (since $x=0$ is a tangent at $(0:1:0)$), and
finally $a_{20}$ must not be zero since otherwise the line $x=0$ would be
contained in $B$ (we, thus, can set $a_{20}=1$). Computing the Jacobian of
the resulting seven equations in the point $(\ell,H)$ (i.e.
for $l_1=l_2=0$) yields the matrix

\[ \left( \begin {array}{cc} 0&0\\-4\,a_{41}&0\\0&0\\0&0
\\0&0\\0&-4\,a_{11}\\0&0\end {array}\right) \]

Thus, $(\ell,H)$ is a singular point in its fibre if and only if
$a_{41}=0$ or $a_{11}=0$, i.e., if and only if $\ell$ touches $B \cap H$
in singular points of $B \cap H$. We only have to compute the ramification
index of the points $(\ell,H)$ which are singular in their fibre. Suppose,
that the line $\ell$ (given by $x=0$) contains the point $(0:0:1)$ which
is assumed to be a singular point of $B \cap H$ and touches $B \cap H$ in
the smooth point $(0:1:0)$. In particular we get the condition $a_{11}=0$.
We get seven polynomials in $l_1$ and $l_2$, the sixth of which has the
form

\[ \left( \cdots\; \mbox{terms containing $l_1$ or $l_2$}\;\cdots - 4a_{22}
      + a_{21}^2 \right) l_2^2.  \]

But $4a_{22} - a_{21}^2$ must not vanish since otherwise $B\cap H$ would
have a cusp in $(0:0:1)$. Thus $(\cdots\; \mbox{terms containing $l_1$ or
$l_2$}\;\cdots - 4a_{22} + a_{21}^2)$ is a unit in
$\complex[l_1,l_2]_{(0,0)}$. Therefore, this equation can be replaced by
$l_2^2$ and the terms containing $l_2^2$ can be cancelled in the other
equations without changing the ideal generated by the equations. Now, the
second equation has the form

\[ \left( \cdots\; \mbox{terms containing $l_1$ or $l_2$}\;\cdots -4\,a_{41}
   \right) l_1.  \]

Since $a_{41}$ must not vanish (otherwise $(0:1:0)$ would be singular on
$B\cap H$), we can replace this equation by $l_1$ and cancel all terms
containing $l_1$ in the other equations. We obtain that the considered
ideal in $\complex[l_1,l_2]_{(0,0)}$ is generated by the two elements
$l_2^2$ and $l_1$. Hence, if $(\ell,H)$ is of the kind that $\ell$
contains the only node of $B\cap H$ then the local ring of $(\ell,H)$ in
its fibre has length two. All other points over $\PPV^3 \setminus \Delta'$
are smooth in their fibres.

For $(\ell,H)\in \PPV^3 \setminus \Delta'$ the fibre $\varphi
^{-1}((\ell,H))$ consists, by Proposition~\ref{Anz-DT}, of 28 points if
$B\cap H$ is smooth, and of 22 points otherwise. By the Remark following
Proposition~\ref{Anz-DT} for a quartic $B\cap H$ with one node, there are
six double tangents through the node.  The Hilbert polynomial of $\varphi
^{-1}((\ell,H)) \in \PPV^2$, hence, is 28 for $B\cap H$ being smooth and
$16+6\cdot 2 = 28$ otherwise. Therefore $\varphi: \YF|_{\PPV^3\setminus
\Delta'} \longrightarrow \PPV^3 \setminus \Delta'$ is flat. From the
length of the local rings computed above we see that the ramification
behaviour is as stated.
\end{proof}

\section{The parameter space of touching conics}\label{sec-tochin-conic}

\subsection{Double tangents and curves in double covers}

Let $B \subset \PP ^3$ be a real quartic surface with exactly 13 nodes,
given by an equation of the form $E_1E_2E_3E_4 -Q^2$, as described in
Proposition~\ref{quartik13}. Let $Z\longrightarrow\PP^3$ be the double
cover which is ramified over $B$. It is constructed as follows (cf.
\cite{bpv}, I.17.). Denote by $p:L\longrightarrow \PP^3$ the total space
of the line bundle $\OO_{\PP^3}(2)$ and by $s\in
H^0(\PP^3,\OO_{\PP^3}(4))$ the section who's zero locus is $B$. Finally,
let $y\in H^0(L,p^*\OO(2))$ be the tautological section. Then $Z\subset L$
is the zero locus of the section $y^2-p^*(s) \in H^0(L,p^*\OO(4))$.

Each of the divisors given by $p^*(E_i)$ splits into two components which
are defined by the sections $y^2 \pm \sqrt{-1}p^*(Q)$ of $H^0(L,p^*\OO(2))$.
Denote these varieties by $S_i^+$ and $S_i^-$, corresponding to $y^2 +
\sqrt{-1}p^*(Q)$ and $y^2 - \sqrt{-1}p^*(Q)$ respectively.

Now, let $H\in \PPV^3$ be a general plane. In particular, let $B_H:=B\cap H$
be a smooth quartic curve. The restriction $Z_H$ of $Z$ to $H$ then is the
smooth double cover of $H\cong\PP^2$ branched along the nonsingular curve
$B_H$. Its canonical bundle is

\[ K_{Z_H}=p^*\OO_H(-3)\otimes p^*\OO_H(2) = p^*\OO_H(-1), \]

which implies that the morphism $p|_H$ is induced by the linear system
$|-K_{Z_H}|$.

By \cite{griffiths-harris} Chapter~4.4, $Z_H$ is isomorphic to the blow-up
of $\PP^2$ in seven points. The 56 $(-1)$-curves are: the seven
exceptional divisors $E^i$ $(i=1,\ldots,7)$, the strict transforms of the
cubics in $\PP^2$ through all seven points with a node in the $i$-th point
$K^i$ $(i=1,\ldots,7)$, the strict transforms of the lines through the
$i$-th and the $j$-th point $G^{ij}$ $(1\le i < j\le 7)$, and the strict
transforms of the conics through all but the $i$-th and the $j$-th point
$C^{ij}$ $(1\le i < j\le 7)$.

The projection $p|_H$ maps the 56 $(-1)$-curves onto the 28 double
tangents of $B_H$ in such a way that each double tangent has exactly two
$(-1)$-curves in its preimage.

The restriction of the $S_i^{\pm}$ to $Z_H$ yields eight curves which are
denoted by $D_i^{\pm}$. These are mapped onto double tangents of $B_H$ and,
hence, are $(-1)$-curves.

\begin{lemma}
$Z_H$ can be realised as the blow-up of $\PP^2$ in such a way that $D_1^+ =
K^7$, $D_1^- = E^7$, $D_2^+=G^{12}$, $D_2^-=C^{12}$, $D_3^+=G^{34}$, $D_3^-
= C^{34}$, $D_4^+=G^{56}$, and $D_4^- = C^{56}$.
\end{lemma}

\begin{proof}
All we have to prove is that the $D_i^{\pm}$ intersect in the correct way,
i.e., $D_i^+ \cdot D_i^- = 2$ for $i=1,\ldots,4$ and $D_i^+\cdot D_j^+ =
D_i^-\cdot D_j^- = 1$ as well as $D_i^+\cdot D_j^- = 0$ for $i\ne j$.

{}From the fact that $p|_H$ is induced by the linear system of the
anticanonical divisor we deduce that $D_i^+ + D_i^- = -K_{Z_H}$ . Hence

\[ 2 = \left(-K_{Z_H}\right)^2 = \left(D_i^+ + D_i^-\right)^2 \]

and consequently $D_i^+ \cdot D_i^- = 2$.

Next observe that the rational function

\[ \frac{y + \sqrt{-1}Q}{(p|_H)^*(E_iE_j)}  \]

corresponds to the principal divisor $D_k^+ + D_l^+ - D_i^- - D_j^-$,
i.e. $[D_k^+ + D_l^+] =[ D_i^- + D_j^-]$ where $\{k,l\} = \{1,\ldots,4\}
\setminus \{i,j\}$. This yields

\[\renewcommand{\arraystretch}{1.5}
   \begin{array}{rcccl}
      \left(D_k^+ + D_l^+\right)^2 & = & \left(D_i^- + D_j^-\right)^2 & = &
         \left(D_k^+ + D_l^+\right)\left(D_i^- + D_j^-\right)\\
      2D_k^+ D_l^+ - 2 & = & 2D_i^- D_j^- -2 & = &
         D_k^+ D_i^- + D_k^+ D_j^- + D_l^+ D_i^- + D_l^+ D_j^- \ge 0
   \end{array}
\]

as the product of different effective divisors is always non-negative.
Hence $D_i^+\cdot D_j^+ = D_i^-\cdot D_j^- = 1$ since two $(-1)$-curves have
intersection product 2 if and only if their sum is an element of
$|-K_{Z_H}|$. But then the sum $D_k^+ D_i^- + D_k^+ D_j^- + D_l^+ D_i^- +
D_l^+ D_j^-$ must vanish and so $D_i^+\cdot D_j^- = 0$.
\end{proof}

Consider now ${\cal Z}:= Z\times_{\PP^3} {\cal H}$ where ${\cal H}
\subset \PP^3 \times \PPV^3$ is the universal (hyper-)plane. Via ${\cal Z}
\rightarrow {\cal H} \rightarrow \PPV^3$, $\cal Z$ is the family of all
surfaces $Z_H$, $H\in \PPV^3$. Let  $\Delta \subset \PPV ^3$ be the set of
all those planes that do not intersect $B$ transversally (i.e., who's
intersection with $B$ is not a smooth quartic curve). Then ${\cal
Z}|_{\PPV^3 \setminus \Delta} \longrightarrow \PPV^3 \setminus \Delta$ is
smooth and proper and therefore, by the Ehresmann--Fibration--Theorem
(cf.  \cite{lamottke}), locally trivial as a fibration of differentiable
manifolds. For any $H_0\in \PPV^3 \setminus \Delta$ the fundamental group
$\pi_1(\PPV^3 \setminus \Delta,H_0)$ acts via monodromy on $H^2({\cal
Z}|_{H_0},\Z) = H^2(Z_{H_0},\Z) = \Pic(Z_{H_0})$ preserving the
intersection pairing. In particular the fundamental group acts via
monodromy on the set of $(-1)$-curves preserving their intersection
behaviour.

On the other hand, this fundamental group $\pi_1(\PPV^3 \setminus
\Delta,H_0)$ acts via monodromy of the finite unramified cover
$\YF|_{\PPV^3\setminus\Delta}\longrightarrow \PPV^3\setminus\Delta$ on the
set of double tangents of $B\cap H_0$. (Recall that \YF{} was defined as
$\YF:=\{ (\ell,H) \in \Grass(2,4)\times \PPV^3\,|\, \ell\in H\,\mbox{and
$\ell$ is double tangent at $B$}\}$.) Obviously, this action is the same as
the action that is induced from the action on $Z_{H_0}$ by mapping a
$(-1)$-curve to its corresponding double tangent.

The properties of the monodromy action on the $(-1)$-curves are to be
examined in the sequel. First, observe that the eight curves $D_i^{\pm}$ keep
fixed under the monodromy action since they are restrictions of the globally
defined $S_i^{\pm}$.

A pair of $(-1)$-curves lying over the same double tangent is mapped to a
pair of $(-1)$ curves that again are mapped onto the same (maybe
different) double tangent. This follows from the fact that the sum
$[C+C']$ of such a pair $(C,C')$ is equal to the anti-canonical class.
which for the $(-1)$-curves is equivalent to $C\cdot C'=2$.

The set of $(-1)$-curves that are different from, say, $D_1^{\pm}$ can be
split in two monodromy--invariant subsets each containing 27 curves: one
subset consisting of all curves $C$ with $C\cdot D_1^- = 1$ and the other
containing the curves $C$ satisfying $C\cdot D_1^- = 0$. I.e., these
subsets are characterised by the intersection number of their elements
with $D_1^- = E^7$ which can take the values $0$ and $1$. As $D_1^-$ is
invariant under the monodromy action this condition is invariant and so
are the subsets. Each pair $(C,C')$ of $(-1)$-curves with $C\cdot C'=2$
(except the pair $(D_1^+,D_1^-)$) has one member in each of the two sets.
Therefore the monodromy action on the $(-1)$ curves is determined by the
action on one of the invariant subsets.

We choose the set of curves that do not intersect $E^7$. This set with its
incidence relations is equivalent to the set of the 27 lines of a smooth
cubic surface: We get the correspondence by blowing down $E^7$. The
$(-1)$-curves not intersecting $E^7$ are mapped onto the $(-1)$-curves in
the blown-down surface which is $\PP^2$ blown-up in six points.

\subsection{The lines in a cubic surface}

Let $S$ be the cubic surface obtained by blowing down $E^7\subset
Z_{H_0}$.  Denote the lines in $S$ by $E^i$ $(i=1,\ldots, 6)$ for the
images of $E^i\subset Z_{H_0}$, $G^{ij}$ $(1\le i<j\le 6)$ for the images
of the corresponding curves $G^{ij}$ in $Z_{H_0}$, and $C^i$ $(i=1,\ldots
6)$ for the images of the curves $C^{i7} \subset Z_{H_0}$.  The monodromy
action on the $(-1)$-curves of $Z_{H_0}$ induces an action of the
fundamental group $\pi_1(\PPV^3\setminus\Delta,H_0)$ on the 27 lines of
$S$ and, hence, induces a morphism of $\pi_1(\PPV^3\setminus\Delta,H_0)$
into the group of symmetries of lines in the cubic surface $S$ (i.e. the
group of permutations of the 27 lines that respect their incidence
relations).

Let $G$ be the image of the fundamental group in this symmetry group.
Clearly, $G$ leaves the three lines $G^{12}$ ($\hateq D_2^+$), $G^{34}$
($\hateq D_3^+$), and $G^{56}$ ($\hateq D_4^+$) invariant. To each of them
there are exactly 10 lines that intersect this line. Furthermore, to each
pair $(C,C')$ of intersecting lines there is exactly one line that
intersects them both. Thus, to each of the three lines $G^{12}$, $G^{34}$,
and $G^{56}$ there is associated a set of eight lines that intersect this
curve and that are different from these three lines. The three sets of
eight lines have to be disjoint since the three lines $G^{12}$, $G^{34}$,
and $G^{56}$ meet each other and so the unique line that intersects two of
them is just the third.

\begin{satz}\label{transitiv}
$G$ acts transitively on the three $G$--invariant sets of curves
associated to the three curves $G^{12}$, $G^{34}$, and $G^{56}$.
\end{satz}

\begin{proof}
The variety \YF{} is naturally fibred over \YO{} -- the parameter space of
double tangents. The fibres are isomorphic to $\PP^1$ since each line in
$\PP^3$ is contained in a pencil of planes. Therefore the irreducible
components of \YF{} are just the preimages of the irreducible components
of \YO{}. The projection $\YF{}\longrightarrow\PPV^3$ is a finite cover
which is \'etale over $\PPV^3\setminus\Delta$ by Proposition~\ref{flach}
and, hence, $\YF|_{\PPV^3 \setminus \Delta}$ is smooth. Therefore its
irreducible components are just its arc--connected components in the
Euclidian topology. By Proposition~\ref{DT-compon}, \YF{} has seven
components. A general fibre $\YF|_H$ over $\PPV^3 \setminus \Delta$
contains 28 points by Proposition~\ref{Anz-DT}, four of them corresponding
to the four double tangents in the planes $E_i=0$. (Recall that the closed
subset in \YO{} of lines in one of the planes $E_i=0$ were recognised to
be irreducible components.) The remaining 24 points belong to double
tangents of the other three components of \YO{}. We claim that any fibre
$\YF|_H$ over $\PPV^3$ must contain at least one point of each
component.

For each component of \YO{} there is a plane in $\PPV^3 \setminus \Delta$
which contains a double tangent of this component: For a double tangent
$\ell$ that does not pass through one of the singular points of $B$, the
pencil of planes containing $\ell$ is not contained in $\Delta \subset
\PPV^3$. The set of double tangents through a singular point $p\in B$ is
one-dimensional as this set is parametrised by the ramification locus of
the projection of $B$ from $p$.  On the other hand, each component of
\YO{} is two-dimensional and, hence, in any component of \YO{} there is an
open set of double tangents $\ell$ not through any of the singular points
of $B$. Now, choose a component of \YO{} and a double tangent $\ell$ of
this component which does not contain singular points of $B$.  Let $H\in
\PPV^3 \setminus \Delta$ be a plane containing $\ell$, i.e.  $(\ell,H) \in
\YF|_{\PPV^3 \setminus \Delta}$. Since $\YF|_{\PPV^3 \setminus \Delta}
\longrightarrow \PPV^3 \setminus \Delta$ is an \'etale cover the component
of \YF{} containing $(\ell,H)$ dominates $\PPV^3$.  Therefore, every plane
in $\PP^3$ contains at least one double tangent of each component of
\YO{}.

As mentioned above, the fundamental group $\pi_1(\PPV^3 \setminus
\Delta,H_0)$ acts on the fibre $\YF|_{H_0}$ via monodromy and, as the
irreducible components of $\YF|_{\PPV^3 \setminus \Delta}$ are just its
connected components, each orbit of the monodromy action contains all
those points of the fibre that belong to the same component.

On the other hand, the monodromy action on $\YF|_{H_0}$ is obtained from
the monodromy action on $(-1)$-curves on $Z|_{H_0}$ by projecting the
$(-1)$-curve onto the corresponding double tangent. So the action on the
27 lines of a cubic surface cannot have more than six orbits. (The seventh
orbit in $\YF|_{H_0}$ consisting of the point corresponding to the double
tangent in the plane $E_1=0$ has no representation among the 27 lines.)
The three sets $\{G^{12}\}$ $\{G^{34}\}$ $\{G^{56}\}$ are $G$-invariant
and consequently the three $G$--invariant sets of lines intersecting
$G^{12}$, $G^{34}$ or $G^{56}$ have to be $G$-orbits.
\end{proof}

\begin{folg}
Let $H$ be a plane in $\PP^3 \setminus \Delta$. There are four components
of \YO{} with exactly one double tangent in $H$ -- namely the four
components parametrising the lines in the planes $E_i=0$ ($i=1, \ldots,
4$). Of each of the remaining three components of \YO{} there are eight
double tangents in $H$.
\end{folg}

\begin{proof}
The 24 double tangents in $H$ which are not contained in one of the four
planes $E_i=0$ correspond to the 24 lines in the cubic surface that are
not fixed under the action of $G$. These 24 lines in the cubic split into
three $G$-orbits, each orbit containing eight of them. The $G$-action on
the lines in the cubic is equivalent to the monodromy action on the double
tangents in $H$. (The double tangent in the plane $E_1=0$ is fixed by the
monodromy.) Hence, the 24 double tangents are spread over the three
components in such a way that each component contains eight of them.
\end{proof}

\begin{satz}\label{generated}
$G$ is generated by elements $g$ with the following properties:
\begin{itemize}
\item $g$ is of order two.
\item $g$ leaves at least 15 of the 27 $(-1)$-curves fixed.
\end{itemize}
In other words: There are at most 6 pairs of lines such that the lines in
these pairs are swapped by $g$.
\end{satz}

\begin{proof}
Let $L$ be a general line in $\PPV ^3$. Without loss of generality one
can assume that $H_0 \in L$. Then, due to \cite{lamottke}
Section~7.4.1, $\pi_1(L\setminus \Delta,H_0) \longrightarrow
\pi_1(\PPV^3\setminus \Delta,H_0)$ is surjective and therefore it is
sufficient to study the monodromy of paths in $L \setminus \Delta$. The
fundamental group $\pi_1(L\setminus \Delta,H_0)$ is generated by the
homotopy classes of paths that loop once counter clockwise round one of
the points in $L \cap \Delta$.

By a sufficiently general choice of $L$, we can achieve that for each $H
\in L\cap \Delta$ the plane quartic $B\cap H$ has exactly one ordinary
node as its only singularity.  For a smooth surface $B$ a generic line in
$\PPV ^3$ will do (cf.  \cite{lamottke} Section~1.6.4). The proof in
\cite{lamottke} works as well in the case of hypersurfaces with isolated
singularities. In the proof one only needs to replace the hypersurface by
the open set of its regular points in all occurrences. Hence if for a
plane $H$ of a generic pencil of planes the quartic curve $B\cap H$ has
more or other singularities than one ordinary node then $H$ must contain a
singular point of $B$ (which is an ordinary node in our case). But the set
of those planes $H$, that contain a singular point of $B$ and for which
$B\cap H$  is not a curve with exactly one ordinary node, has at least
codimension two in $\PPV ^3$. So we can choose $L$ in such a way that $L$
does not intersect this codimension-2-subset in $\PPV^3$.

Now, by Proposition~\ref{flach}, the fibre of \YF{} over any $H \in L\cap
\Delta$ contains exactly six ramification points and the ramification
index in each of them is two. The monodromy of a loop round $H$ can only
interchange two sheets of $\YF|_L$ which meet in one of the ramification
points over $H$. Therefore, the monodromy of this loop can swap at most
six couples of points in the fibre of \YF{} over $H_0$. Using the
correspondence between the monodromy action on $\YF|_{H_0}$ and the
monodromy action on the lines in the cubic $S$, this proves the
proposition.
\end{proof}

The group of symmetries of the 27 lines in a cubic surface (i.e. the group
of permutations that respect the incidence relations among the lines) has
been intensively studied. (Cf.~\cite{henderson}, \cite{segre}, \cite{manin}
-- to mention only a few.) The following theorem holds (\cite{manin} Ch.~IV
Theorem~1.9):

\begin{theorem}
The group of symmetries of the 27 lines on a cubic surface is isomorphic to
the Weyl--Group \Esechs.\proved
\end{theorem}

In particular, the group $G$ is a subgroup of \Esechs{}. Using the
Propositions~\ref{transitiv} and \ref{generated} we will be able to
determine $G \subset \Esechs$ as a subgroup of \Esechs{}. For this
purpose, we first determine the elements of \Esechs{} that admit the
properties required in Proposition~\ref{generated}. In \cite{manin}
(Ch.~IV \S~9) as well as in \cite{swinn} a complete list of the conjugacy
classes of \Esechs{} can be found. Moreover, to each conjugacy class the
action of its elements on the 27 lines is described. It turns out that the
only conjugacy class who's elements act on the 27 lines as postulated in
Proposition~\ref{generated} is the class which is denoted by $C_{16}$ in
\cite{manin} and \cite{swinn}.  This class contains exactly 36 elements.

It is a classical result (cf. e.g. \cite{segre}) that the group of
symmetries of the 27 lines in a cubic surface is generated by elements
which swap the lines in one of Schl\"afli's 36 ``double six''.  A ``double
six'' consists of a pair of sextuples of lines such that the lines in each
sextuple are mutually skew and each line of one of the tuples intersects
exactly five lines of the other tuple. Associating to each line of one
tuple the unique line of the other tuple that is skew to this line yields
a one-to-one correspondence between the two sextuples of a double six. We
identify a double six with the element of \Esechs{} that exchanges the
lines of the sextupels in such a way that each line is swapped with the
unique line of the other tuple which is skew to this line. This
transformation keeps the other 15 lines fixed. An example of a double six
is the pair $[(E^1,\ldots,E^6), (C^1,\ldots,C^6)]$.

Obviously, a permutations of the 27 lines corresponding to a double six
satisfies the conditions of Proposition~\ref{generated}, i.e. is in the
conjugacy class $C_{16}$. As this class contains exactly 36 elements
(cf.~\cite{manin} or \cite{swinn}) and as there are exactly 36 double
sixes, $C_{16}$ consists just of these double six transformations.
Consequently, the group $G$ has its generators among these special
transformations.

Next, we will give a list of all double six transformations, that leave
the three lines $G^{12}$, $G^{34}$, and $G^{56}$ invariant. For this
purpose we will use the notation of \cite{manin} and give the
transformations in terms of reflections of a root system. The Picard group
$\Pic(S)$ of a smooth cubic surface $S$ is the free abelian group with
generators $[H]$ -- the pull-back of $\OO_{\PP^2}(1)$ under the blow-up
morphism $S \longrightarrow \PP^2$ -- and the six classes $[-E^i]$ (where
$E^i$ $(i=1,\ldots,6)$ are the exceptional divisors). An element
$a[H]-b_1[E^1]-\cdots - b_6[E^6] \in \Pic(S)$\label{picard} will be
denoted by $(a;b_1,\ldots, b_6)$ in the sequel. The intersection pairing
on the Picard group is given by

\[ (a;b_1,\ldots,b_6)\cdot (a';b'_1,\ldots,b'_6) =
   aa' - \sum_{i=1}^6 b_ib'_i.       \]

Denote by $\omega = -(3;1,1,1,1,1,1)$ the canonical line bundle of the
cubic surface. Let $\omega^\perp \subset \Pic(S) \otimes \R$ be the
orthogonal complement of $\omega$ with respect to the scalar product
induced by the intersection pairing. Note that on $\omega^\perp$ the
intersection pairing is negative definite (cf.~\cite{manin}
Proposition~IV.3.3). There is a root system of type \Esechs{} in
$\omega^{\perp}$ who's  reflections are in one-to-one correspondence with
the double six transformations.  The roots of this root system are the
following:

\begin{itemize}
\item $[E^i] - [E^j] \in \Pic(S) \subset \Pic(S) \otimes \R$ $(1\le i,j
      \le 6, i\ne j)$ (30 roots).
\item $\pm([H]-[E^i]-[E^j]-[E^k])$, $(1\le i < j < k \le 6)$ (40 roots).
\item $\pm (2;1,1,1,1,1,1)$ (2 roots).
\end{itemize}

For a root $x$ the corresponding reflection $s_x$ is given by

\[   s_x(v) = v + (v,x)\,x  \]

(where the scalar product is the one given by the intersection pairing on
$\Pic(S)$). The restriction of any reflection to $\Pic(S) \subset \Pic(S)
\otimes\R$ induces an endomorphism of $\Pic(S)$ that respects the
intersection pairing. In particular it induces a transformation of the 27
lines which is a double six transformation.

Now it is easy to find the double six transformations which fix the three
lines $G^{12}$, $G^{34}$, and $G^{56}$. These correspond to the roots

\begin{equation}\label{roots}
   \renewcommand{\arraystretch}{1.5}
   \begin{array}{lll}

      x_1=  (2;1,1,1,1,1,1), & x_2 =  (1;1,0,1,0,0,1), &
         x_3 =  (0;-1,1,0,0,0,0), \\
      x_4 =  (1;1,0,0,1,1,0),&  x_5 =  (1;1,0,0,1,0,1), &
         x_6 =  (1;0,1,1,0,1,0),\\
      x_7 =  (0;0,0,-1,1,0,0), & x_8 =  (1;0,1,1,0,0,1), &
         x_9 =  (1;0,1,0,1,1,0), \\
      x_{10} =  (1;0,1,0,1,0,1),& x_{11} =  (0;0,0,0,0,-1,1),&
          x_{12} =  (1;1,0,1,0,1,0).\\
   \end{array}
\end{equation}

The corresponding double sixes are:

\begin{equation}\label{action}
\renewcommand{\arraystretch}{2}
\newenvironment{ar}{\renewcommand{\arraystretch}{1.0}%
                    \begin{array}{*{6}{c@{\;}}}}{\end{array}}
\begin{array}{ll}

   x_1 \hateq \left( \begin{ar}
      E^1 & E^2 & E^3 & E^4 & E^5 & E^6 \\
      C^1 & C^2 & C^3 & C^4 & C^5 & C^6
   \end{ar} \right),&
   x_2 \hateq \left( \begin{ar}
      E^1 & E^3 & E^6 & G^{45} & G^{25} & G^{24} \\
      G^{36} & G^{16} & G^{13} & C^2 & C^4 & C^5
   \end{ar} \right),\\[1ex]

   x_3 \hateq \left( \begin{ar}
      E^1 & C^1 & G^{23} & G^{24} & G^{25} & G^{26}\\
      E^2 & C^2 & G^{13} & G^{14} & G^{15} & G^{16}
   \end{ar} \right),&
   x_4 \hateq \left( \begin{ar}
      E^1 & E^4 & E^5 & G^{36} & G^{26} & G^{23} \\
      G^{45} & G^{15} & G^{14} & C^2 & C^3 & C^6
   \end{ar} \right),\\[1ex]

   x_5 \hateq \left( \begin{ar}
      E^1 & E^4 & E^6 & G^{35} & G^{25} & G^{23} \\
      G^{46} & G^{16} & G^{14} & C^2 & C^3 & C^5
   \end{ar} \right),&
   x_6 \hateq \left( \begin{ar}
      E^2 & E^3 & E^5 & G^{46} & G^{16} & G^{14} \\
      G^{35} & G^{25} & G^{23} & C^1 & C^4 & C^6
   \end{ar} \right),\\[1ex]

   x_7 \hateq \left( \begin{ar}
      E^3 & C^3 & G^{14} & G^{24} & G^{45} & G^{46}\\
      E^4 & C^4 & G^{13} & G^{23} & G^{35} & G^{36}
   \end{ar} \right),&
   x_8 \hateq \left( \begin{ar}
      E^2 & E^3 & E^6 & G^{45} & G^{15} & G^{14} \\
      G^{36} & G^{26} & G^{23} & C^1 & C^4 & C^5
   \end{ar} \right),\\[1ex]

   x_9 \hateq \left( \begin{ar}
      E^2 & E^4 & E^5 & G^{36} & G^{16} & G^{13} \\
      G^{45} & G^{25} & G^{24} & C^1 & C^3 & C^6
   \end{ar} \right),&
   x_{10} \hateq \left( \begin{ar}
      E^2 & E^4 & E^6 & G^{35} & G^{15} & G^{13} \\
      G^{46} & G^{26} & G^{24} & C^1 & C^3 & C^5
   \end{ar} \right),\\[1ex]

   x_{11} \hateq \left( \begin{ar}
      E^5 & C^5 & G^{16} & G^{26} & G^{36} & G^{46}\\
      E^6 & C^6 & G^{15} & G^{25} & G^{35} & G^{45}
   \end{ar} \right),&
   x_{12} \hateq \left( \begin{ar}
      E^1 & E^3 & E^5 & G^{46} & G^{26} & G^{24} \\
      G^{35} & G^{15} & G^{13} & C^2 & C^4 & C^6
   \end{ar} \right).
\end{array}
\end{equation}

The 12 roots (\ref{roots}) (together with their negatives) form a root
system. The group $G'$ which is generated by the double six
transformations corresponding to these roots is the Weyl group to this
root system. Obviously $G\subset G'$ is a subgroup. To determine $G'$
observe that the four roots $x_3,x_7,x_{11}$ and $x_{12}$ form a basis of
the root system (\ref{roots}), i.e. any other root in the system is a
linear combination of these four roots and the coefficients are either all
non-negative or all non-positive. $G'$ is, thus, uniquely determined by
the corresponding Dynkin diagram:

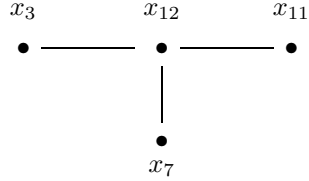
\begin{center}
\unitlength=0.7pt
\begin{picture}(145.00,92.00)(115.00,615.00)
\put(125.00,680.00){\line(1,0){50.00}}
\put(200.00,680.00){\line(1,0){50.00}}
\put(190.00,670.00){\line(0,-1){30.00}}
\put(115.00,680.00){\circle*{6}}
\put(190.00,680.00){\circle*{6}}
\put(260.00,680.00){\circle*{6}}
\put(190.00,630.00){\circle*{6}}
\put(115.00,700.00){\makebox(0,0)[cc]{$x_3$}}
\put(190.00,700.00){\makebox(0,0)[cc]{$x_{12}$}}
\put(260.00,700.00){\makebox(0,0)[cc]{$x_{11}$}}
\put(190.00,615.00){\makebox(0,0)[cc]{$x_7$}}
\end{picture}
\end{center}

This diagram clearly is the Dynkin diagram of the group $D_4$. The group
$D_4$ is isomorphic to the semi-direct product of the permutation group
$S_4$ with $(\Z_2)^3$.

We are, now, going to show that $G$ must be the whole group $G'$. For this
purpose consider the curve $G^{12}$ which is fixed under the group
action.  As already mentioned, there are ten lines in the cubic $S$ that
intersect $G^{12}$.  Since for any two intersecting lines in a cubic
surface there exists a unique line in this cubic which intersects them
both the ten lines which intersect $G^{12}$ are grouped into five pairs of
intersecting lines. These are $p_0:=(G^{34},G^{56})$ (which is fixed under
the monodromy action) and $p_1:=(G^{46},G^{35})$, $p_2:=(G^{36},G^{45})$,
$p_3:=(C^1,E^2)$, $p_4:=(C^2,E^1)$.

As the group $G'$ acts on the 27 lines preserving their incidence
relations and keeping the line $G^{12}$ fixed, each of the above pairs of
intersecting lines must be moved to such a pair by the group action. Using
the explicit description (\ref{action}) of the action of $G'$ one observes
that any of the transformations corresponding to the $x_i$ interchanges
two of the four pairs $p_1,\ldots,p_4$ -- either preserving or reversing
the order of the lines in the pairs. As $G\subset G'$ is to act
transitively on the lines of the pairs $p_1,\ldots,p_4$ it must, in
particular, act transitively on these four pairs. But a transitive
subgroup of $S_4$ that is generated by transpositions must be the group
$S_4$ itself. (Remember that $G$ is to be generated by a subset of the
transformations corresponding to the $x_i$.)

Let $x_{i_1}$, $x_{i_2}$, and $x_{i_3}$ correspond to elements of $G$ that
act by swapping the pairs $p_1\leftrightarrow p_2$, $p_2\leftrightarrow
p_3$, and $p_3\leftrightarrow p_4$ respectively. The subgroup of $G$
generated by these elements is isomorphic to $S_4$. The action of any of
its elements does not map one line in a pair $p_i$ onto the other line of
this pair. Therefore -- as $G$ is to act transitively on the eight lines
-- there is an $x_{i_4}'$ such that the corresponding transformation is in
$G$ and is not in the subgroup generated by the $x_{i_j}$ $(j=1,2,3)$. Let
$x_{i_4}$ be the transformation that acts on the pairs in the same manner
as $x_{i_4}'$ does and which is in the subgroup generated by the $x_{i_j}$
$(j=1,2,3)$.  Then the composition of $x_{i_4}$ and $x'_{i_4}$ swaps the
lines in the two pairs that are interchanged by the action of $x_{i_4}$
and $x'_{i_4}$.

By conjugating $x_{i_4}'$ with the elements of
$\left<x_{i_1},x_{i_2},x_{i_3}\right> \subset G$ on can interchange the
lines in any two of the four pairs.  Hence $G$ contains a subgroup
(isomorphic to $(\Z_2)^3$) who's elements interchange the lines in an even
number of pairs. This proves that $G$ contains a subgroup which is the
semi direct product of $S_4$ and $(\Z_2)^3$ and therefore $G=G'$.

Now, that we have a detailed knowledge of the action on the 27 lines (or
on the 56 $(-1)$-curves on $Z_{H_0}$) induced by the monodromy, we will
examine the induced action on the pairs of (different) lines formed out of
these lines. Denote by \OA{}, \OB{}, \OC{} the three orbits consisting of
the eight lines intersecting the lines $G^{12}$, $G^{34}$, $G^{56}$
respectively. Then there are the following monodromy invariant subsets in
the set of pairs of lines:

\begin{itemize}
\item The 9 sets $\{G_{12}\}\times\OA$, $\{G_{12}\}\times\OB$,
      $\{G_{12}\}\times\OC$, $\{G_{34}\}\times\OA$ etc.
\item The 3 sets $\{(G_{12},G_{34})\}$, $\{(G_{12},G_{56})\}$, and
      $\{(G_{34},G_{56})\}$.
\item The 3 sets $(\OA\OA)$ (which means the set of pairs consisting of
      two lines of \OA{}), $(\OB\OB)$, $(\OC\OC)$.
\item The 3 sets $(\OA\OB)$ (which means the set of pairs consisting of
      one line of \OA{} and one line of \OB{}), $(\OA\OC)$, $(\OB\OC)$.
\end{itemize}

The sets of the first two items are obviously orbits under the monodromy
action, whereas the sets of the last two items will turn out to be
composite of two orbits.

One orbit in the set $(\OA\OA)$ consists of the four pairs
$p_1:=(G^{46},G^{35})$, $p_2:=(G^{36},G^{45})$, $p_3:=(C^1,E^2)$,
$p_4:=(C^2,E^1)$ which are the four pairs which consist of intersecting
lines. We claim that the remaining 24 elements in $(\OA{}\OA)$ form an
orbit of the monodromy action.

Let $(\ell_1,\ell_2)\in \OA$ be a pair of lines which do not intersect. We
will calculate the number of elements in the orbit of this pair by finding
its stabiliser. Embed $S_4\subset G$ as the subgroup who's elements
permute the pairs $p_1,\ldots,p_4$ but leave the order of the lines in the
pairs unchanged. If $(\Z_2)^3 \subset G$ is the subgroup of elements that
leave the four pairs $p_1,\ldots,p_4$ fixed but swaps the lines in an even
number of these pairs then $G$ is the semi-direct product of $S_4$ and
$(\Z_2)^3$ (maybe in a different presentation as above). Let the lines
$\ell_1$ and $\ell_2$ belong to the pairs $p_{i_1}$ and $p_{i_2}$
respectively. Then the stabiliser of $(\ell_1, \ell_2)$ consists of those
elements of $G$ which leave the pairs $p_{i_1}$ and $p_{i_2}$ unchanged
(i.e. those which interchange the two other pairs with or without changing
the order of the lines in the pairs) together with those elements that
interchange $p_{i_1}$ and $p_{i_2}$ in such a way that $\ell_1$ is mapped
to $\ell_2$. The stabiliser is, therefore, $(S_2\times S_2) \ltimes \Z_2
\subset S_4 \ltimes (\Z_2)^3 \cong G$. The orbit of $(\ell_1,\ell_2)$,
hence must contain $24=192/8$ elements. Applying the same argument to
\OB{} and \OC{} yields: \OA{}, \OB{}, and \OC{} each split into two orbits
-- one with four and one with 24 elements.

In the same manner we will attack the sets $(\OA\OB)$, $(\OA\OC)$, and
$(\OB\OC)$ -- e.g. the set $(\OA\OB)$ (the two other sets being treated in
an analogous way). Consider the line $E^1\in\OA$. The stabiliser in $G$ of
$E^1$ consists of exactly those elements which act only on the three pairs
different from $(C^2,E^1)$. This subgroup is isomorphic to
$S_3 \ltimes (\Z_2)^2$. It is generated by the elements corresponding to the
double six transformations $x_6,\ldots,x_{11}$. Using the explicit
description (\ref{action}) it is easy to check that the set \OB{} has two
orbits under the action of this subgroup -- namely
$\{E^3,E^4,G^{25},G^{26}\}$ and $\{C^3,C^4,G^{15},G^{16}\}$. (Note that
the second orbit consists of the set of lines in \OB{} that intersect
$E^1$.) Thus, any pair $(\ell_1,\ell_2) \in (\OA\OB)$ can be moved into a
pair $(E^1,\ell)$ via monodromy action and this pair can be moved into
each of the pairs $(E^1,\ell')$ with $\ell' \in\{E^3,E^4,G^{25},G^{26}\}$
or $\ell' \in\{C^3,C^4,G^{15},G^{16}\}$ depending on in which set $\ell$
is contained. So $(\OA\OB)$ can contain at most two orbits. On the other
hand, a pair $(\ell_1,\ell_2)$ of intersecting lines cannot be moved into
a pair of non-intersecting lines and vice versa since the monodromy action
respects the incidence relations.  Hence, $(\OA\OB)$ must split in at
least two orbits -- one containing pairs of intersecting lines and the
other containing pairs of non-intersecting lines. This proves that
$(\OA\OB)$ (as well as $(\OA\OC)$ and $(\OB\OC)$) are the composite of two
orbits of equal cardinality. Summing up, we have proved the following
proposition.

\begin{satz}\label{DT-Paare-monod}
The set of pairs formed out of the 27 lines has the following orbits under
the action of $G$ induced by the monodromy action on the lines.
\begin{enumerate}
\item Each of the three subsets $(\OA\OA)$, $(\OB\OB)$, $(\OC\OC)$ contains
      two orbits -- one with four and one with 24 pairs.
\item Each of the subsets $(\OA\OB)$, $(\OA\OC)$, $(\OB\OC)$ contains two
      orbits with 32 pairs -- one orbit with pairs of intersecting lines
      and one orbit with pairs of non-intersecting lines.
\item The 9 sets $\{G_{12}\}\times\OA$, $\{G_{12}\}\times\OB$,
      $\{G_{12}\}\times\OC$, $\{G_{34}\}\times\OA$ etc. with 8
      pairs each and
\item the 3 sets $\{(G_{12},G_{34})\}$, $\{(G_{12},G_{56})\}$, and
      $\{(G_{34},G_{56})\}$ -- are orbits.\proved
\end{enumerate}
\end{satz}

\subsection{One-parameter-families of touching conics and linear systems}

We, now, want to use our knowledge on pairs of lines of a cubic surface
and their monodromy to examine the irreducible components of the parameter
space of touching conics. For this purpose return to the the double cover
$Z$ of $\PP^3$ branched along our quartic $B$.

Let $H\subset\PP^3$ be an arbitrary plane that intersects $B$
transversally and denote by $Z_H$ the restriction of $Z$ to $H$ and by
$\pi:Z_H \longrightarrow H$ the induced morphism. (Recall that $Z_H$ is
isomorphic to the blow-up of $\PP^2$ in seven points in general position.)
We will establish a connection between one-parameter-families of touching
conics of $B\cap H$ and certain linear systems in $Z_H$. Thereby the
reducible elements of the linear systems will correspond to the reducible
touching conics in the one-parameter-families.

Let $C_1$ and $C_2$ be two $(-1)$-curves that have different images under
$\pi$ (i.e. $\pi(C_1)\ne\pi(C_2)$) which means that they are $(-1)$-curves
over different double tangents of $B\cap H$. This is equivalent to
$[C_1+C_2] \ne -K_{Z_H}$ where $K_{Z_H}$ denotes the canonical class of
$Z_H$. Let $C'_1$ and $C'_2$ be the $(-1)$-curves defined by $[C_i + C'_i]
= -K_{Z_H}$. Then these curves intersect as follows:

\[ C_1\cdot C_2 = C'_1\cdot C'_2 = 1 - C_1\cdot C'_2 = 1 - C'_1\cdot C_2. \]

This follows from

\[ 1 = C_1 \cdot (-K_{Z_H}) = C_1\cdot\left( C_2+C'_2\right)
     = C_1\cdot C_2 + C_1\cdot C'_2 \]

and analogous identities. Therefore, by eventually interchanging $C_1$ and
$C'_1$, one can achieve that $C_1\cdot C_2 = 1$ (the corresponding double
tangents keeping unchanged).

\begin{satz}
If $C_1$ and $C_2$ are chosen as above with $C_1\cdot C_2 =1$ then the
linear system $|C_1+C_2|$ is one-dimensional. Its generic element is a
smooth rational curve that by the projection $\pi:Z_H \longrightarrow H$
is mapped onto a touching conic.
\end{satz}

\begin{proof}
For any of the 56 $(-1)$-curves $C$ of $Z_H$ consider the exact sequence

\[  0\longrightarrow \OO_{Z_H}\longrightarrow \OO_{Z_H}(C)
     \longrightarrow \OO_{C}(C) \longrightarrow 0.\]

Noting that $C$ is a smooth rational curve with $\OO_{C}(C)=
\OO_{C}(C\cdot C) =\OO_{C}(-1)$ and that $Z_H$ is a smooth rational
surface so that $0=h^1(\OO_{Z_H}) = h^1(\OO_{Z_H}(C))$ we get
$h^0(\OO_{Z_H}) = h^0(\OO_{Z_H}(C)) = 1$. Now, the exact sequence

\[  0\longrightarrow \OO_{Z_H}(C_2) \longrightarrow \OO_{Z_H}(C_1 + C_2)
     \longrightarrow \OO_{C_1}(C_1 + C_2) \longrightarrow 0  \]

and the fact that $\OO_{C_1}(C_1 + C_2) = \OO_{C_1}$ (since $C_1\cdot(C_1+
C_2) = 0$) yield

\[ h^0(\OO_{Z_H}(C_1 + C_2)) = h^0(\OO_{Z_H}(C_2)) + h^0(\OO_{C_1}) = 2 \]

and, hence, $\dim |C_1+C_2|=1$.

$|C_1 +C_2|$ cannot have a fixed component. Since $C_1$ and $C_2$ are
irreducible this fixed component would have to be one of these two curves
and then the linear system would only contain the divisor $C_1+C_2$ in
contradiction to the dimension of the system being 1. Hence, as $(C_1
+C_2)^2=0$, the system cannot have base points at all.

By Bertini's Theorem the generic element of $|C_1 +C_2|$ is smooth away from
the base locus. As the base locus is empty the generic element of the
system is smooth everywhere. Let $C\in |C_1 +C_2|$ be general and $D_1,
\ldots, D_n$ be its irreducible components. The $D_i$ cannot intersect each
other for any intersection point would be a singular point of $C$. Thus
$D_i\cdot D_j = 0$ for $i\ne j$ and consequently

\[ 0 = \left( C_1 +C_2\right)^2 = \left( \sum_{i=1}^n D_i\right)^2 =
       \sum_{i=1}^n D_i^2.  \]

$|C_1 +C_2|$ has no fixed components -- therefore

\[  0\le D_i\cdot (C_1 +C_2) = D_i^2 \]

and hence $D_i^2 = 0$ for all $i$. From the adjunction formula we get

\[ g_i := \mbox{genus}(D_i) = \frac{K_{Z_H}\cdot D_i}{2} +1. \]

On the other hand $K_{Z_H} \cdot D_i < 0$ since $-K_{Z_H}$ is ample. From
$g_i \ge 0$ we then get $K_{Z_H} \cdot D_i = -2$ for all~$i$. But

\[  -K_{Z_H} \cdot \left(\sum_{i=1}^n D_i\right) =
    -K_{Z_H}\cdot \left( C_1 +C_2\right) = 2 \]

and hence $n=1$, $C=D_1$, and $\mbox{genus}(C)=0$, i.e. $C$ is a smooth
rational curve.

Now, let $R\subset Z_H$ be the ramification divisor of the map
$Z_{H}\longrightarrow H$: $\OO_{Z_H}(R) = \pi^*(\OO_H(2)) = -2K_{Z_H}$.
Hence

\[ (C_1 +C_2)\cdot [R] = 2\,(C_1 + C_2)\cdot(-K_{Z_H}) = 4  \]

and by projection formula

\[ 4\,[pt] = \pi_*\left(\left(C_1+C_2\right)\cdot [R]\right)
         = \pi_*\left(\left(C_1+C_2\right) \cdot
               \pi^*\left(2\,[\,l\,]\right)\right)
         = \pi_*\left(C_1+C_2\right)\cdot \left(2\,[\,l\,]\right) \]

where $[pt]$ denotes the class of a point and $[\,l\,]$ the class of a
line in $H$. Therefore $\pi_*(C_1+C_2)\in |2\,[\,l\,]|$. Let $C\in
|C_1+C_2|$ be a general element. If $\pi(C)$ were a line then $C$ were
contained in the preimage of a line. As $[C_1+C_2] \ne K_{Z_H}$ there would
exist an effective $C'$ such that $[C]+[C'] = -K_{Z_H}$. But then
$C'\cdot(-K_{Z_H}) = (-K_{Z_H})^2 - C\cdot(-K_{Z_H}) = 0$ which is
impossible as $-K_{Z_H}$ is ample. So the the image of $C$ under $\pi$
must be a smooth conic and $\pi|_C$ is of degree one onto $\pi(C)$.

Consider now the preimage $\pi^{-1}(\pi(C))$ in $Z_H$. As $\pi$ is a
double cover and $\pi|_C$ is only of degree one, $\pi^{-1}(\pi(C))$ must
contain other components than $C$ or $\pi(C)$ must be contained in $B\cap
H$. The latter is not possible since $B\cap H$ was supposed to be a smooth
quartic curve. Obviously, $\pi^{-1}(\pi(C))$ is an element of
$|-2\,K_{Z_H}|$ since $\pi(C)$ is an element of $|\OO_{\PP^2}(2)|$ and
$\pi:Z_H\longrightarrow \PP^2\cong H$ is induced by the anticanonical
linear system. Therefore, the sum of the other components of
$\pi^{-1}(\pi(C))$ must be an element of $|-2\,K_{Z_H} - (C_1 +C_2)| =
|C'_1 +C'_2|$. ($C'_i$ was defined to be the $(-1)$-curve in $Z_H$ such
that $C_i + C'_i = -K_{Z_H}$.)

Let $C'\in |C'_1 + C'_2|$ be the divisor which is complementary to $C$ in
$\pi^{-1}(\pi(C))$. Note that $(C'_1 +C'_2)$, as well as $(C_1+C_2)$, is
the sum of two $(-1)$-curves with intersection $C'_1 \cdot C'_2 =1$ and
consequently the above arguments equally apply to $(C'_1 +C'_2)$. In
particular, a general element of $|C'_1 +C'_2|$ is a smooth rational curve
which by $\pi$ is mapped onto a smooth conic in $H$. So if $C \in |C_1
+C_2|$ is sufficiently general then $C'$, as well, is a smooth rational
curve that is mapped onto a smooth conic. Therefore, $\pi^{-1}(\pi(C))$
splits into two components each of which is a smooth rational curve.

Now, one shows just like in the proof of Proposition~\ref{geraden-in-DS}
that the preimage in a double cover of a conic in $H$ splits into two
components if and only if it has even intersection with the ramification
locus $B\cap H \subset H$. Therefore $C$ is mapped onto a touching conic
and the proposition is proved.
\end{proof}

By the above proposition $\pi$ induces a morphism $|C_1+C_2|\cong \PP^1
\longrightarrow \PP^5$ where $\PP^5$ is the parameter space of conics in
$H$. This morphism is necessaryly injective as the preimage of $\pi(C)$
consists of $C$ and an element of the linear system $|-K_{Z_H} - [C]|$
which is different from $|C_1 +C_2|$. There is an open subset in $|C_1
+C_2|$ which is mapped into the closed subset of touching conics in
$\PP^5$. Therefore any element of $|C_1 +C_2|$ is mapped onto a (maybe
reducible) touching conic. Moreover, an element of $|C_1 +C_2|$ is mapped
to a reducible conic if and only if it is the sum of two $(-1)$-curves.

We, thus, have constructed a correspondence between the linear systems
$|C_1 + C_2|$ (with $(-1)$-curves $C_i$ satisfying $C_1 \cdot C_2 =1$) and
one-parameter-families of conics in $H$ touching $B\cap H$: For each one
parameter family there exist exactly two of these linear systems that are
mapped to this family.

Let again ${\cal Z}:= Z\times_{\PP^3} {\cal H}$ where ${\cal H} \subset
\PP^3 \times \PPV^3$ is the universal (hyper-)plane. We have seen that for
any $H_0\in \PPV^3 \setminus \Delta$ the fundamental group
$\pi_1(\PPV^3 \setminus \Delta,H_0)$ acts via monodromy of ${\cal
Z}|_{\PPV^3\setminus\Delta} \longrightarrow \PPV^3\setminus\Delta$ on
$\Pic(Z_{H_0})$ preserving the intersection pairing. In particular, the
fundamental group acts on the above liner systems.

Denote by $\X'\subset P$ ($P$ -- the parameter space of all conics in
$\PP^3$ as constructed in Section~\ref{psc}) the closed subscheme of
conics that have only even intersection with $B$ and let $\X\subset\X'$ be
the union of all irreducible components that do not entirely consist of
double lines. From $P$ \X{} inherits a morphism to $\PPV^3$. By
Proposition~\ref{Anz-epf} the fibre of \X{} over any $H\in
\PPV^3\setminus\Delta$ consists of 63 disjoint smooth conic curves in the
fibre of $P$ over $H \in \PPV^3$ which is isomorphic to $\PP^5$. The
fundamental group $\pi_1(\PPV^3 \setminus \Delta,H_0)$ acts on the set of
the 63 one parameter families in the fibre $\X_{H_0}$ in a natural way by
monodromy:  Any path $\gamma:[0,1]\longrightarrow \PPV^3 \setminus\Delta$
can be lifted (in a non-unique way) to a path in
$\X|_{\PPV^3\setminus\Delta}$. Though the lift is not unique, the
connected component of the fibre $\X_{\gamma(t)}$ in which the lifted path
is contained is well determined.

Obviously, the correspondence between the linear systems in $Z_{H_0}$ and
the one-parameter-families is compatible with the monodromy action. On the
other hand, the orbits of the monodromy in $\X_{H_0}$ are in bijection
with the connected components of $\X|_{\PPV^3\setminus\Delta}$. But
$\X|_{\PPV^3\setminus\Delta}$ is smooth (since it is flat over
$\PPV^3\setminus\Delta$ and has smooth fibres) and therefore the connected
components of $\X|_{\PPV^3\setminus\Delta}$ in the Euclidian topology are
its irreducible components.

We are mainly interested in the connected components of
$\X|_{\PPV^3\setminus\Delta}$. One way of computing the monodromy orbits
would be to list all linear systems $|C_1+C_2|$ with $(-1)$-curves $C_i$
satisfying $C_1\cdot C_2 =1$ and then determining the monodromy orbits of
these linear systems using our knowledge on the monodromy of
$(-1)$-curves. The monodromy of $\X|_{\PPV^3\setminus\Delta}$ is then
easily calculated. This approach is a bit cumbersome. So we modify this
method using our knowledge about one-parameter-families of touching
conics.

{\sloppy
First, we consider the linear system $|G^{12}+G^{34}| =
|(2;1,1,1,1,0,0,0)|$ (the elements of $\Pic(Z_{H_0})$ are denoted similarly
as on page~\pageref{picard}, i.e. $(a;b_1,\ldots,b_7) \in \Pic(Z_{H_0})$
denotes the element $a[H] - b_1[E^1] - \cdots - b_7[E^7]$).
As the two $(-1)$-curves $G^{12}$ and $G^{34}$ are monodromy invariant the
linear system must keep fixed under the monodromy action, as well.  This
system may be presented as the sum of two $(-1)$-curves in five further
ways:}

\[ \renewcommand{\arraystretch}{1.5}
\begin{array}{rcccccc}
   (2;1,1,1,1,0,0,0) & = & C^{56}+E^7 & = & C^{57}+E^6 & = & C^{67}+E^5\\
                     & = & G^{13}+G^{24} & = & G^{14}+G^{23}.
\end{array} \]

The two pairs $(G^{12},G^{34})$ and $(C^{56},E^7)$ correspond to the pairs
of double tangents $(e_2,e_3)$ and $(e_4,e_1)$ where $e_i$ denotes the
double tangent given by the equation $E_i=0$ in $H_0$ ($E_i$
being the linear forms appearing in the equation $E_1E_2E_3E_4 -Q^2$ of
$B$).  The remaining four pairs are those of the set $(\OA\OA)$ when
identified with the corresponding lines in a cubic surface.  The
corresponding one-parameter-family of touching conics, thus, contains four
pairs of double tangents all eight double tangents being in the same
component of \YO. Denote this component of \YO{} by \OC{}.

Analogously (considering the systems $|G^{12}+G^{56}|$ and $|G^{34} +
G^{56}|$ respectively) the one-parameter-family of touching conics
containing the pairs $(e_2,e_4)$ and $(e_1,e_3)$ contains four pairs
formed of double tangents of the component \OB{} of \YO{} and the family
containing the pairs $(e_3,e_4)$ and $(e_1,e_2)$ contains four pairs
formed of double tangents of the component \OA{}.  These three
one-parameter-families are the ones which are obvious from the special
form of the equation of $B$: $E_1E_2E_3E_4 -Q^2$ is of the form $UW-V^2$
by letting $V=Q$ and letting $U$ be one of $E_1E_2$, $E_1E_3$ or $E_1E_4$.
So in any plane $H_0$ we get three one-parameter-families~(\ref{epf})
which we will call the ``obvious'' families.

In particular, each of these obvious one-parameter-families contains two
reducible conics consisting of the four double tangents $e_1,\ldots,e_4$
and four reducible elements formed of the eight double tangents in $H_0$
that belong to the same component \OA{}, \OB{} or \OC{} of \YO{}.

Denote, for simplicity the eight double tangents in $H_0$ of the component
\OA{} by $a_1,\ldots,a_8$ and the double tangents of \OB{} and \OC{} by
$b_1,\ldots,b_8$ and $c_1,\ldots,c_8$ respectively.  We will examine how
the pairs formed out of these double tangents can be distributed to the
one-parameter-families. By Lemma~\ref{epf-zerf-el} each one-parameter-family of
touching conics in $H_0$ contains exactly six reducible conics. From the
above arguments we already know the reducible elements of three families
(given in terms of pairs of double tangents):

\[
\renewcommand{\arraystretch}{1.5}
\begin{array}{*{3}{r@{,\hspace{\arraycolsep}}}*{3}{l@{\hspace{\arraycolsep}}}}
   e_1e_2& e_3e_4& a_1a_2& a_3a_4,&a_5a_6,&a_7a_8 \\
   e_1e_3& e_2e_4& b_1b_2&\ldots \\
   e_1e_4& e_2e_3& c_1c_2&\ldots \\
\end{array}
\]

In particular, the equation of the plane quartic $B\cap H_0$ can be
written in the form $e_1e_2\cdot a_1a_2 -V^2$ where $V$ is a quadratic
form. (Lines and the corresponding linear forms are denoted by the same
letter.) By writing this equation as $e_1a_1\cdot e_2a_2 -V^2$ we see that
the couples $e_1a_1$ and $e_2a_2$ belong to the same one-parameter-family.
By Lemma~\ref{DT-in-EPF}, in this family no further couple of double
tangents is contained that has one of the double tangents $e_i$ or $a_i$
as an element.  Couples of the form $b_ib_j$ or $c_ic_j$ also must not
occur in this family.  If, for instance, the couple $b_1b_3$ were in one
group together with $e_1a_2$ then the equation of $B\cap H_0$ could be
written as $e_1a_2 \cdot b_1b_3 - V'^2 = e_1b_1\cdot a_1b_3 - V'^2$ and
thus $e_1b_1$ and $a_1b_3$ would belong to the same family. But this is
impossible by the above argument.  Hence, only the couples $b_kc_l$ may
occur in families together with couples $e_ia_j$. There are exactly 16
families each of which contains two couples $e_ia_j$.  On the other hand,
there are 64 couples $b_ic_j$ which is just the number of couples needed
to complete these 16 families.

By the same argument, the couples $a_ib_j$ belong to families which
contain two couples of the form $e_kc_l$ and four couples of the form
$a_ib_j$ and, analogously, the couples $a_ic_j$ spread over families
containing two couples $e_kb_j$ and four couples $a_ic_j$.

As we have just seen couples $b_ic_j$ and $b_kc_l$ have to occur in one
family. Thus, there is a family containing both $b_ib_k$ and $c_jb_l$.
Analogously, there are couples $a_ia_j$ and $b_kb_l$ as well as $a_ia_j$
and $c_kc_l$ in one group. There are $3\cdot({8 \choose 2} -4) = 72$ pairs
of the form $a_ia_j$, $b_ib_j$ and $c_ic_j$ which do not occur in the
``obvious'' families. These spread over the remaining 12
one-parameter-families.\myfootnote{In fact, using Lemma~\ref{DT-in-EPF},
one could determine which pairs of double tangents pertain to which
one-parameter-family. Those considerations can be found in \cite{salmon}.}

In the consequence of these considerations we can determine the
connected components of $X|_{\PPV^3 \setminus\Delta}$. Let $\YY\subset
X|_{\PPV^3 \setminus\Delta}$ be the closed subscheme which parametrises
the reducible touching conics. \YY{} is naturally isomorphic to the open
subset in the symmetric product of \YF{} with itself:

\[  \YF\symm_{\PPV^3\setminus \Delta} \YF \setminus \mbox{Diag}
    \;\tilde{\longrightarrow} \;\YY  \]

by simply associating to each pair of complanar double tangents the
corresponding reducible touching conic. As any one-parameter-family in any
fibre of $X|_{\PPV^3 \setminus\Delta}$ over $\PPV^3\setminus\Delta$
contains reducible conics and as the monodromy action on the set of
one-parameter-families in the fibre over $H_0 \in \PPV^3\setminus\Delta$
is independent of the chosen lift of the path one can choose the lift of
any path $\gamma\subset \PPV^3\setminus\Delta$ in such a way that the
lifted path is contained in \YY{}.

But the connected components of \YY{} are already determined: The following
proposition is just a corollary of Proposition~\ref{DT-Paare-monod}.

\begin{satz}\label{YY-comp}
The connected components of $\YY{}\longrightarrow \PPV^3\setminus
\Delta$ are the following:

\begin{enumerate}
\item \label{orb24}Three components with 24 points in each fibre over
      $\PPV^3 \setminus \Delta$ -- each component corresponding to one
      orbit (the one with 24 elements) in the sets $(\OA\OA)$, $(\OB\OB)$,
      and $(\OC\OC)$ of pairs of lines in a cubic.
\item \label{orb4}Three components with four points in each fibre -- each
      component corresponding to the other orbit in the sets $(\OA\OA)$,
      $(\OB\OB)$, and $(\OC\OC)$ of pairs of lines in a cubic.
\item Six components with 32 points in each fibre corresponding to the
      orbits in the sets $(\OA\OB)$, $(\OA\OC)$, and $(\OB\OC)$.
\item \label{orb8}12 components with eight points in each fibre: For
      each $i=1,\ldots 4$ and each component \OA{}, \OB{} or \OC{}
      of \YO{} there is an irreducible component of \YY{}. The
      corresponding reducible conics in a plane $H$ consist of the
      double tangent $e_i$ given by the equation $E_i=0$ in $H$ and
      one double tangent of the chosen component of \YO{}.
\item Six components with just one point in each fibre, namely the six
      pairs of double tangents $e_ie_j$.
\end{enumerate}
\end{satz}

Consequently, $X|_{\PPV^3 \setminus\Delta}$ has the following connected
components:

\begin{itemize}
\item The three ``obvious'' families in the fibre of $X|_{\PPV^3
      \setminus\Delta}$ over $H_0 \in\PPV^3 \setminus\Delta$ are
      invariant under monodromy as they contain reducible conics
      consisting of the double tangents $e_ie_j$ which are fixed
      under monodromy. Each of the three families contains two of
      them. The other four reducible conics in each family are
      necessarily the four pairs of double tangents of one of the
      orbits in item~\ref{orb4}.\ in the above Proposition.
\item There are six connected components with eight
      one-parameter-families in any fibre of $X|_{\PPV^3
      \setminus\Delta}$ over $\PPV^3 \setminus\Delta$: All
      one-parameter-families containing reducible conics of the type
      $e_ia_j$, $e_ib_j$ or $e_ic_j$ pertain to one of these
      components. As we have seen, any one-parameter-family that
      contains those pairs of double tangent contains two of them
      and four pairs of the form $a_ib_j$, $a_ic_j$, or $b_ic_j$. By
      Lemma~\ref{DT-in-EPF} and the above discussion, the two pairs
      $e_ia_j$ and $e_ka_l$ are in the same one-parameter-family
      only if $i\ne k$ and $j\ne l$ (analogously for $e_ib_j$ and
      $e_ic_j$). Therefore each component of \YY{} in
      item~\ref{orb8}.\ of Proposition~\ref{YY-comp} intersects a
      one-parameter-family in a fibre over $\PPV^3 \setminus\Delta$
      in at most one point.  Consequently, the monodromy orbit of
      such a one-parameter-family in the fibre of $X|_{\PPV^3
      \setminus\Delta}$ over $H_0$ consists of exactly eight
      families.
\item The remaining pairs of double tangents are those of
      item~\ref{orb24}.\ of Proposition~\ref{YY-comp}. The families
      containing these reducible conics belong to the same connected
      component: We have seen that there is a one-parameter-families
      that contains a pair $a_ia_j$ as well as a pair $b_kb_l$ and a
      family that contains some $a_ia_j$ together with a $c_kc_l$.
      So, one can connect the first family with any family
      containing a pair $a_\bullet a_\bullet$ or a pair $b_\bullet
      b_\bullet$ by a path lying in \YY{} and, analogously, the second
      family can be connected with any family containing a pair $c_\bullet
      c_\bullet$.  Hence the twelve one-parameter-families belong to the
      same connected component.
\end{itemize}

As the connected components of $X|_{\PPV^3 \setminus\Delta}$ are just its
irreducible components the following Theorem is proved by the above
discussion.

\begin{theorem}
$X|_{\PPV^3 \setminus\Delta}$ has 10 irreducible components -- namely
\begin{itemize}
\item three components each with one one-parameter-family in every fibre
      over $\PPV^3$,
\item six components with eight families in every fibre, and
\item one component with twelve families in each fibre over $\PPV^3
      \setminus\Delta$.\proved
\end{itemize}
\end{theorem}

The three types of irreducible components differ by the type of reducible
conics that they contain. In each one-parameter-family there are reducible
conics that contain two, one or none double tangent $e_i$ respectively for
the three types.

{\bf Remark:} The components of $X$ are not entirely determined by the
above theorem. There are at least four components which are contained in
$X|_\Delta$ (i.e. over $\Delta \subset \PPV^3$) namely the four sets
consisting of all conics in the planes $E_i=0$. All of them have to be
touching conics since $B\cap{\{E_i=0\}}$ is a non-reduced conic. But these
four components are, conjecturally, all components which are contained in
$X|_\Delta$.

\vspace{\fill}
{\small
Ingo Hadan\\
Institut f\"ur Reine Mathematik\\
Humboldt-Universit\"at zu Berlin\\
Ziegelstra\ss{}e 13a\\
10099 Berlin\\
e-mail: hadan@mathematik.hu-berlin.de
}

\end{document}